\documentclass[MSNbibl,number,dvips]{arxstspdf}
\usepackage{flushend}
\usepackage{stfloats}
\usepackage{graphicx}

\volume{26}
\issue{1}
\pubyear{2011}
\firstpage{84}
\lastpage{101}
\doi{10.1214/10-STS347}

\makeatletter

\def\@bmisc[#1]{%
  \get@battribute{unstr}%
  \common@pub@types%
  \let\bauthor\bbl@bauthor%
  \let\bhowpublished\@firstofone%
  \def\borganization##1{{\bauthor@style ##1}}%
}

\newproclaim{answer}{Answer}
\newtheorem{theorem}{Theorem}
\newtheorem{proposition}{Proposition}

\renewcommand{\emptyset}{\varnothing}
\renewcommand{\epsilon}{\varepsilon}

\newcommand{\st}{\mid}
\newcommand{\givn}{\vert}

\newcommand{\dd}{\mathrm{d}}            

\newcommand{\CCC}{\mathcal{C}}          
\newcommand{\FFF}{\mathcal{F}}          
\newcommand{\GGG}{\mathcal{G}}          

\newcommand{\III}{\mathbb{I}}          

\newcommand{\bbbr}{\mathbb{R}}                  

\newcommand{\Prob}{\mathbf{P}}
\newcommand{\Qrob}{\mathbf{Q}}

\newcommand{\Expect}{\mathbf{E}}
\newcommand{\var}{\mathbf{var}}

\makeatother

\begin{document}
\begin{frontmatter}

\title{Test Martingales, Bayes Factors and~$p$-Values}
\runtitle{Martingales and $p$-Values}

\begin{aug}
\author{\fnms{Glenn} \snm{Shafer}\ead[label=e1,text=gshafer@rbsmail.rutgers.edu]{gshafer@rbsmail.rutgers.edu}},
\author{\fnms{Alexander} \snm{Shen}\ead[label=e2]{alexander.shen@lif.univ-mrs.fr}},
\author{\fnms{Nikolai} \snm{Vereshchagin}\ead[label=e3]{ver@mech.math.msu.su}}
\and
\author{\fnms{Vladimir} \snm{Vovk}\corref{}\ead[label=e4]{vovk@cs.rhul.ac.uk}}
\runauthor{Shafer, Shen, Vereshchagin and Vovk}

\affiliation{Rutgers Business School, University of London, CNRS,
Moscow State University and University of London}

\address{Glenn Shafer is Professor,
  Rutgers Business School, Newark, New  Jersey, USA and Department of~Com\-puter Science,
  Royal Holloway, University of London,
  Egham, Surrey, United Kingdom \printead{e1}.
  Alexander Shen is Directeur de Recherche,
  Laboratoire d'Informatique Fondamentale, CNRS,
  Marseille, France \printead{e2}.
  Nikolai Vereshchagin is Professor,
  Department of Mathematical Logic and the Theory of Algorithms,
  Faculty of Mechanics and Mathematics,
  Moscow State University, Moscow~119899, Russia \printead{e3}.
  Vladimir Vovk is Professor, Department of Computer Science,
  Royal Holloway, University of London,
  Egham, Surrey, United Kingdom \printead{e4}.}

\end{aug}

\begin{abstract}
 A nonnegative martingale with initial value equal to one
  measures evidence against a probabilistic hypothesis.
  The inverse of its value at some stopping time can be interpreted as a Bayes factor.
  If we exaggerate the evidence
  by considering the largest value attained so far by such a martingale,
  the exaggeration will be limited,
  and there are systematic ways to eliminate it.
  The inverse of the exaggerated value at some stopping time
  can be interpreted as a $p$-value.
  We give a simple characterization
  of all increasing functions that eliminate the exaggeration.
\end{abstract}

\begin{keyword}
\kwd{Bayes factors}
    \kwd{evidence}
    \kwd{hypothesis testing}
    \kwd{martingales}
    \kwd{$p$-values}.
\end{keyword}

\end{frontmatter}

\section{Introduction}
\label{sec:introduction}

\begin{figure*}

\includegraphics{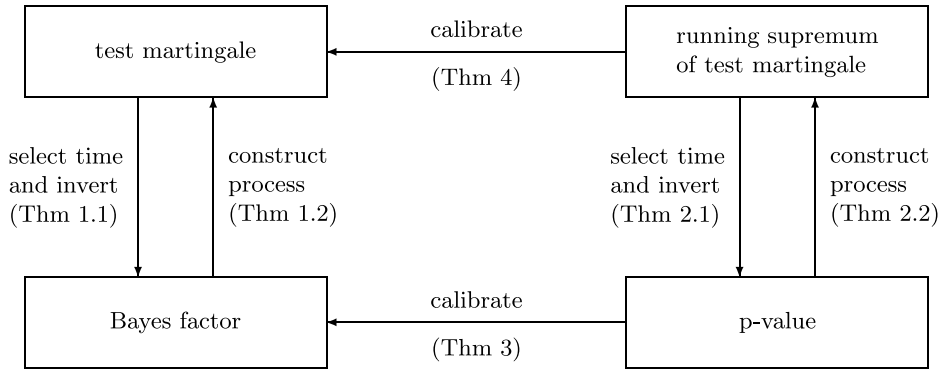}

\caption{The relationship between a Bayes factor and a $p$-value
  can be thought of as a snapshot of the dynamic relationship
  between a nonnegative martingale $(X_t)$ with initial value $1$
  and the process $(X^*_t)$ that tracks its supremum.
  The snapshot could be taken at any time,
  but in our theorems we consider the final
  values of the martingale
  and its supremum process.}\label{fig:analogy}
\end{figure*}

Nonnegative martingales with initial value $1$,\break Bayes factors and
$p$-values can all be regarded as measures of evidence against a
probabilistic hypothesis (i.e., a~simple statistical hypothesis). In
this article we review the well-known relationship between Bayes
factors and nonnegative martingales and the less well-known
relationship between $p$-values and the suprema of nonnegative
martingales. Figure \ref{fig:analogy} provides a visual frame for the
relationships we\break discuss.

Consider a random process $(X_t)$ that initially has the value one and
is a nonnegative martingale under a~probabilistic hypothesis $P$ (the
time $t$ may be discrete or continuous). We call such a martingale a
\textit{test martingale}.
One statistical interpretation of the values of a test martingale is
that they measure the changing evidence against $P$. The value $X_t$ is
the number of dollars a gambler has at time $t$ if he begins with $\$1$
and follows a certain strategy for betting at the rates given by~$P$;
the nonnegativity of the martingale means that this strategy never
risks a cumulative loss exceeding the $\$1$ with which it began. If
$X_t$ is very large, the gambler has made a lot of money betting
against $P$, and this makes $P$ look doubtful. But then $X_u$ for some
later time $u$ may be lower and make $P$ look better.

The notion of a test martingale $(X_t)$ is related to the notion of a
Bayes factor, which is more familiar to statisticians. A Bayes factor
measures the degree to which a fixed body of evidence supports $P$
relative to a particular alternative hypothesis $Q$; a very small value
can be interpreted as discrediting $P$. If $(X_t)$ is a test
martingale, then for any fixed time $t$, $1/X_t$ is a Bayes factor. We
can also say, more generally, that the value $1/X_{\tau}$ for any
stopping time $\tau$ is a Bayes factor. This is represented by the
downward arrow on the left in Figure~\ref{fig:analogy}.

Suppose we exaggerate the evidence against~$P$ by considering not the
current value $X_t$ but the greatest value so far:
\[
  X^*_t := \sup_{s\le t} X_s.
\]
A high $X^*_t$ is not as impressive as a high $X_t$, but how should we
understand the difference? Here are two complementary answers:

\begin{answer}[(Downward arrow on the right in
Figure~\ref{fig:analogy})]\label{answer1}
  Although $(X^*_t)$ is usually not a martingale,
  the final value $X^*_{\infty}:=\sup_s X_s$
  still has a property associated with hypothesis testing:
  for every $\delta\in[0,1]$, $1/X^*_{\infty}$
  has probability no more than $\delta$ of being $\delta$ or less.
  For any $t$, $X^*_t$, because it is less than or equal to $X^*_{\infty}$,
  has the same property.
  In this sense, $1/X^*_{\infty}$ and $1/X^*_t$ are $p$-values (perhaps conservative).
 \end{answer}

\begin{answer}[(Leftward arrow at the top of Figure~\ref{fig:analogy})]\label{answer2}
  As we will show, there are systematic ways of shrinking
  $X^*_t$  (\textit{calibrating} it, as we shall say)
  to eliminate the exaggeration.
  There exist, that is to say, functions~$f$
  such that $\lim_{x\to\infty}f(x)=\infty$
  and $f(X^*_t)$ is an unexaggerated measure of evidence\break against $P$,
  in as much as there exists a test martingale $(Y_t)$
  always satisfying $Y_t \ge f(X^*_t)$ for all $t$.
\end{answer}

Answer \ref{answer2} will appeal most to readers familiar with the algorithmic
theory of randomness, where the idea of treating a martingale as a
dynamic measure of evidence is well established (see, e.g.,
\cite{livitanyi2008}, Section~4.5.7). Answer \ref{answer1} may be more interesting
to readers familiar with mathematical statistics, where the static
notions of a~Bayes factor and a $p$-value are often compared.

For the sake of conceptual completeness, we note that Answer \ref{answer1} has a
converse. For any random variable $p$ that has probability $\delta$ of
being $\delta$ or less for every $\delta\in[0,1]$, there exists a test
martingale $(X_t)$ such that $p=1/X^*_{\infty}$. This converse is
represented by the upward arrow on the right of our figure. It may be
of limited practical interest, because the time scale for $(X_t)$ may
be artificial.

Parallel to the fact that we can shrink the running supremum of a test
martingale to obtain an unexaggerated test martingale is the fact that
we can inflate a $p$-value to obtain an unexaggerated Bayes factor.  This
is the leftward arrow at the bottom of Figure \ref{fig:analogy}. It was previously
discussed in \cite{vovk1993} and \cite{sellkeetal2001}.

These relationships are probably all known in one form or another to
many people. But they have received less attention than they deserve,
probably because the full picture emerges only when we bring together
ideas from algorithmic randomness and mathematical statistics. Readers
who are not familiar with both fields may find the historical
discussion in Section \ref{sec:history} helpful.

Although our theorems are not deep, we state and prove them using the
full formalism of modern probability theory.  Readers more comfortable
with the conventions and notation of mathematical statistics may want
to turn first to Section \ref{sec:example}, in which we apply these
results to testing whether a coin is fair.

The theorems depicted in Figure \ref{fig:analogy} are proven in
Sections \ref{sec:math}--\ref{sec:reduce}. Section \ref{sec:math} is
devoted to mathematical preliminaries; in particular, it introduces the
concept of a~test martingale and the wider and, in general, more
conservative concept of a test supermartingale. Section~\ref{sec:bayes}
reviews the relationship between test supermartingales and Bayes
factors, while Section \ref{sec:pvalues} explains the relationship
between the su\-prema of test supermartingales and $p$-values.
Section \ref{sec:calibrate} explains how $p$-values can be inflated so
that they are not exaggerated relative to Bayes factors, and
Section \ref{sec:reduce} explains how the maximal value attained so far
by a test supermartingale can be similarly shrunk so that it is not
exaggerated relative to the current value of a test supermartingale.

There are two appendices.  Appendix \ref{app:continuous} explains why
test supermartingales are more efficient tools than test martingales in
the case of continuous time.  Appendix \ref{app:calculations} carries
out some calculations that are used in Section~\ref{sec:example}.

\section{Some History}\label{sec:history}

Jean Ville introduced martingales into probability theory in his 1939
thesis \cite{ville1939}. Ville considered only test martingales and
emphasized their betting interpretation. As we have explained, a test
martingale under~$P$ is the capital process for a betting strategy that
starts with a unit capital and bets at rates given by $P$, risking only
the capital with which it begins. Such a strategy is an obvious way to
test $P$: you refute the quality of $P$'s probabilities by making money
against them.

As Ville pointed out, the event that a test martingale tends to
infinity has probability zero, and for every event of probability zero,
there is a test martingale that tends to infinity if the event happens.
Thus, the classical idea that a probabilistic theory predicts events to
which it gives probability equal (or nearly equal) to one can be
expressed by saying that it predicts that test martingales will not
become infinite (or very large). Ville's idea was popularized after
World War II by Per Martin-L\"{o}f
\cite{martinlof1966b,martinlof1969} and subsequently developed by
Claus-Peter Schnorr in the 1970s \cite{schnorr1971} and A. P. Dawid in
the 1980s \cite{dawid1984}. For details about the role of martingales
in algorithmic randomness from von Mises to Schnorr, see
\cite{bienvenuetal2009}. For historical perspective on the paradoxical
behavior of martingales when they are not required to be nonnegative
(or at least bounded below), see \cite{brubruchung2009}.

Ville's idea of a martingale was taken up as a~technical tool in
probability mathematics by Joseph Doob in the 1940s \cite{locker2009},
and it subsequently became important as a technical tool in
mathematical statistics, especially in sequential analysis and time
series \cite{lai2009} and in survival analysis \cite{aalenetal2009}.
Mathematical statistics has been slow, however, to take up the idea of
a martingale as a dynamic measure of evidence. Instead, statisticians
emphasize a static concept of hypothesis testing.

Most literature on statistical testing remains in the static and
all-or-nothing (reject or accept) framework established by Jerzy Neyman
and Egon Pearson in 1933 \cite{neymanpearson1933}. Neyman and Pearson
emphasized that when using an observation $y$ to test $P$ with respect
to an alternative hypothesis $Q$, it is optimal to reject $P$ for
values of $y$ for which the likelihood ratio $P(y)/Q(y)$ is smallest
or, equivalently, for which the reciprocal likelihood ratio $Q(y)/P(y)$
is largest. [Here $P(y)$ and $Q(y)$ represent either probabilities
assigned to $y$ by the two hypotheses or, more generally, probability
densities relative to a common reference measure.] If the observation
$y$ is a vector, say, $y_1, \ldots, y_t$, where $t$ continues to grow,
then the reciprocal likelihood ratio $Q(y_1, \ldots, y_t)/P(y_1,\ldots,
y_t)$ is a discrete-time martingale under $P$, but mathematical
statisticians did not propose to interpret it directly. In the
sequential analysis invented by Abraham Wald and George A. Barnard in
the 1940s, the goal still is to define an all-or-nothing Neyman--Pearson
test satisfying certain optimality conditions, although the reciprocal
likelihood ratio plays an important role [when testing $P$ against~$Q$,
this goal is attained by a rule that rejects $P$ when  $Q(y_1,\ldots,\break
y_t)/P(y_1,\ldots, y_t)$ becomes large enough and accepts $P$ when
$Q(y_1,\ldots, y_t)/P(y_1,\ldots, y_t)$ becomes small\break enough].

The increasing importance of Bayesian philosophy and practice starting
in the 1960s has made the likelihood ratio $P(y)/Q(y)$ even more
important. This ratio is now often called the Bayes factor for $P$
against~$Q$, because by Bayes's theorem, we obtain the ratio of $P$'s
posterior probability to $Q$'s posterior probability by multiplying the
ratio of their prior probabilities by this factor
\cite{kassraftery1995}.


The notion of a $p$-value developed informally in statistics. From Jacob
Bernoulli onward, everyone who applied probability theory to
statistical data agreed that one should fix a threshold (later called a
\textit{significance level}) for probabilities, below which a probability
would be small enough to justify the rejection of a~hypothesis. But
because different people might fix this threshold differently, it was
natural, in empirical work, to report the smallest threshold for which
the hypothesis would still have been rejected, and British
statisticians (e.g., Karl Pearson in 1900 \cite{pearson1900} and
R.~A.~Fisher in 1925 \cite{fisher1925}) sometimes called this
borderline probability ``the value of $P$.'' Later, this became
``$P$-value'' or ``$p$-value'' \cite{aldrich2009}.

After the work of Neyman and Pearson, which emphasized the
probabilities of error associated with significance levels chosen in
advance, mathematical statisticians often criticized applied
statisticians for merely reporting $p$-values, as if a small $p$-value were
a~measure of evidence, speaking for itself without reference to a
particular significance level. This disdain for $p$-values has been
adopted and amplified by modern Bayesians, who have pointed to cases
where \mbox{$p$-values} diverge widely from Bayes factors and hence are very
misleading from a Bayesian point of view
\cite{sellkeetal2001,wagenmakers2007}.

\section{Mathematical Preliminaries}\label{sec:math}

In this section we define martingales, Bayes factors and $p$-values. All
three notions have two versions: a~narrow version that requires an
equality and a wider version that relaxes this equality to an
inequality and is considered conservative because the goal represented
by the equality in the narrow version may be more than attained; the
conservative versions are often technically more useful. The
conservative version of a martingale is a supermartingale. As for Bayes
factors and $p$-values, their main definitions will be conservative, but
we will also define narrow versions.

Recall that a \textit{probability space} is a triplet
$(\Omega,\FFF,\Prob)$, where $\Omega$ is a set, $\FFF$ is a
$\sigma$-algebra on $\Omega$ and $\Prob$ is a~probability measure on
$\FFF$. A \textit{random variable} $X$ is a real-valued $\FFF$-measurable
function on $\Omega$; we allow random variables to take values
$\pm\infty$. We use the notation $\Expect(X)$ for the integral of $X$
with respect to $\Prob$ and $\Expect(X\givn\GGG)$ for the conditional
expectation of $X$ given a~$\sigma$-algebra $\GGG\subseteq\FFF$; this
notation is used only when $X$ is integrable [i.e., when
$\Expect(X^+)<\infty$ and $\Expect(X^-)<\infty$; in particular,
$\Prob\{X=\infty\}=\Prob\{X=-\infty\}=0$]. A \textit{random process} is a
family $(X_t)$ of random variables $X_t$; the index $t$ is interpreted
as time. We are mainly interested in discrete time (say,
$t=0,1,2,\ldots$),
but our results (Theorems \ref{thm:bayes}--\ref{thm:reduce}) will also
apply to continuous time (say, $t\in[0,\infty)$).

\subsection{Martingales and Supermartingales}
\label{subsec:supermartingales}

The time scale for a martingale or supermartingale is formalized by a
filtration. In some cases, it is convenient to specify this filtration
when introducing the martingale or supermartingale; in others, it is
convenient to specify the martingale or supermartingale and derive an
appropriate filtration from it. So there are two standard definitions
of martingales and supermartingales in a probability space. We will use
them both:
\begin{longlist}
\item[(1)]
  $(X_t,\FFF_t)$, where $t$ ranges over an ordered set
  ($\{0,1,\ldots\}$ or $[0,\infty)$ in this article),
  is a \textit{supermartingale}
  if $(\FFF_t)$ is a filtration
  (i.e., an indexed set of sub-$\sigma$-algebras of $\FFF$
  such that $\FFF_s\subseteq\FFF_t$ whenever $s<t$),
  $(X_t)$ is a random process adapted with respect to $(\FFF_t)$
  (i.e., each $X_t$ is $\FFF_t$-measurable),
  each $X_t$ is integrable, and
  \[ 
    \Expect(X_t\givn\FFF_s)
    \le
    X_s\quad\mbox{a.s.}
  \]
  when $s<t$.
  A supermartingale is a \textit{martingale} if,
  for all $t$ and $s<t$,
  \begin{equation}\label{eq:martingale}
    \Expect(X_t\givn\FFF_s)
    =X_s\quad\mbox{a.s.}
  \end{equation}
\item[(2)]
  A random process $(X_t)$ is a \textit{supermartingale}
  (resp. \textit{martingale})
  if $(X_t,\FFF_t)$ is a supermartingale (resp. martingale),
  where $\FFF_t$ is the $\sigma$-algebra generated by $X_s$, $s\le t$.
\end{longlist}
For both definitions, the class of supermartingales contains that of
martingales.

In the case of continuous time we will always assume that the paths of
$(X_t)$ are right-continuous almost surely (they will then
automatically have left limits almost surely; see, e.g.,
\cite{dellacheriemeyer1982}, VI.3(2)). We will also assume that the
filtration $(\FFF_t)$ in $(X_t,\FFF_t)$ satisfies the \textit{usual
conditions}, namely, that each $\sigma$-algebra $\FFF_t$ contains all
subsets of all $E\in\FFF$ satisfying $\Prob(E)=0$ (in particular, the
probability space is complete) and that $(\FFF_t)$ is
\textit{right-continuous}, in that, at each time $t$,
$\FFF_t=\FFF_{t+}:=\bigcap_{s>t}\FFF_s$. If the original filtration
$(\FFF_t)$ does not satisfy the usual conditions (this will often be
the case when $\FFF_t$ is the $\sigma$-algebra generated by $X_s$,
$s\le t$), we can redefine $\FFF$ as the $\Prob$-completion
$\FFF^{\Prob}$ of $\FFF$ and redefine $\FFF_t$ as
$\FFF^{\Prob}_{t+}:=\bigcap_{s>t}\FFF^{\Prob}_s$, where $\FFF^{\Prob}_s$
is the $\sigma$-algebra generated by $\FFF_s$ and the sets
$E\in\FFF^{\Prob}$ satisfying $\Prob(E)=0$; $(X_t,\FFF_t)$ will remain
a (super)martingale by \cite{dellacheriemeyer1982}, VI.3(1).

We are particularly interested in \textit{test supermartingales}, defined
as supermartingales that are nonnegative ($X_t \ge 0$ for all $t$) and
satisfy $\Expect(X_0)\le1$, and \textit{test martingales}, defined as
martingales that are nonnegative and satisfy $\Expect(X_0)=1$.
Earlier, we defined test martingales as those having initial value 1;
this can be reconciled with the new definition by setting $X_t:=1$ for
$t<0$.
A well-known fact about test supermartingales, first proven for
discrete time and test martingales by Ville, is that
%
\begin{equation}\label{eq:maximal}
  \Prob\{X^*_{\infty} \ge c \} \le 1/c
\end{equation}
for every $c\ge 1$ (\cite{ville1939}, page 100;
\cite{dellacheriemeyer1982}, VI.1). We will call this the \textit{maximal
inequality}. This inequality shows that $X_t$ can take the value
$\infty$ only with probability zero.

\subsection{Bayes Factors}

A nonnegative measurable function $B\dvtx\Omega\to[0,\infty]$
is called a \textit{Bayes factor for $\Prob$} if $\int(1/B)\,\dd\Prob\le1$;
we will usually omit ``for $\Prob.$'' A Bayes factor $B$ is said to be
\textit{precise} if $\int(1/B)\,\dd\Prob=1$.

In order to relate this definition to the notion of Bayes factor
discussed informally in Sections \ref{sec:introduction}
and~\ref{sec:history}, we note first that whenever $\Qrob$ is a
probability measure on $(\Omega,\FFF)$, the Radon--Nikodym derivative
$\dd\Qrob/\dd\Prob$
will satisfy $\int(\dd\Qrob/\dd\Prob)\,\dd\Prob\le1$, with equality if
$\Qrob$ is absolutely continuous with respect to $\Prob$. Therefore,
$B=1/(\dd\Qrob/\dd\Prob)$ will be a Bayes factor for $\Prob$. The Bayes
factor $B$ will be precise if $\Qrob$ is absolutely continuous with
respect to $\Prob$; in this case $B$ will be a version of the
Radon--Nikodym derivative $\dd\Prob/\dd\Qrob$.

Conversely, whenever a nonnegative measurable function $B$ satisfies
$\int(1/B)\,\dd\Prob\le1$, we can construct a probability measure $\Qrob$
that has $1/B$ as its Radon--Nikodym derivative with respect to $\Prob$.
We first construct a measure $\Qrob_0$ by setting $\Qrob_0(A) := \int_A
(1/B)\, \dd \Prob$ for all $A\in\FFF$, and then obtain $\Qrob$ by adding
to $\Qrob_0$ a~measure that puts the missing mass $1-\Qrob_0(\Omega)$
(which can be $0$) on a set $E$ (this can be empty or a single point)
to which $\Prob$ assigns probability zero. (If $\Prob$ assigns positive
probability to every element of $\Omega$, we can add a new point to
$\Omega$.) The function $B$ will be a version of the Radon--Nikodym
derivative $\dd\Prob/\dd\Qrob$ if we redefine it by setting
$B(\omega):=0$ for $\omega\in E$ [remember that $\Prob(E)=0$].

\subsection{$p$-Values}

In order to relate $p$-values to supermartingales, we introduce a new
concept, that of a $p$-test. A \textit{$p$-test} is a~measurable function
$p\dvtx\Omega\to[0,1]$ such that
%
\begin{equation}\label{eq:p}
  \Prob\{\omega\st p(\omega)\le\delta\}
  \le
  \delta
\end{equation}
for all $\delta\in[0,1]$. We say that $p$ is a \textit{precise $p$-test} if
%
\begin{equation}\label{eq:precise}
  \Prob\{\omega\st p(\omega)\le\delta\}
  =
  \delta
\end{equation}
for all $\delta\in[0,1]$.

  It is consistent with established usage to call the values of a $p$-test \textit{$p$-values},
  at least if the $p$-test is precise.
  One usually starts from a measurable function $T\dvtx\Omega\to\bbbr$
  (the \textit{test statistic})
  and sets $p(\omega):=\Prob\{\omega'\st T(\omega')\ge T(\omega)\}$;
  it is clear that a function $p$ defined in this way,
  and any majorant of such a $p$, will satisfy~(\ref{eq:p}).
  If the distribution of $T$ is continuous,
  $p$ will also satisfy (\ref{eq:precise}).
  If not, we can treat the ties $T(\omega')=T(\omega)$ more carefully
  and set
  \begin{eqnarray*} 
    p(\omega)
    &:=&
    \Prob\{\omega'\st T(\omega')>T(\omega)\}\\
    &&{}+
    \xi\Prob\{\omega'\st T(\omega')=T(\omega)\},
  \end{eqnarray*}
  where $\xi$ is chosen randomly from the uniform distribution on $[0,1]$;
  in this way we will always obtain a~function satisfying (\ref{eq:precise})
  (where $\Prob$ now refers to the overall probability
  encompassing generation of $\xi$).

\section{Supermartingales and Bayes Factors}
\label{sec:bayes}

When $(X_t,\FFF_t)$ is a test supermartingale, $1/X_t$ is a~Bayes
factor for any value of $t$. It is also true that $1/X_{\infty}$,
$X_{\infty}$ being the supermartingale's limiting value, is a Bayes
factor. Part 1 of the following theorem is a~precise statement of the
latter assertion; the former assertion follows from the fact that we
can stop the supermartingale at any time $t$.

Part 2 of Theorem \ref{thm:bayes} states that we can construct a~test martingale
whose limiting value is reciprocal to a given precise Bayes factor. We
include this result for mathematical completeness rather than because
of its practical importance; the construction involves arbitrarily
introducing a filtration, which need not correspond to any time scale
with practical meaning. In its statement, we use $\FFF_{\infty}$ to
denote the $\sigma$-algebra generated by $\bigcup_{t}\FFF_t$.

\begin{theorem}\label{thm:bayes}
  (1)
    If $(X_t,\FFF_t)$ is a test supermar\-tingale,
    then $X_{\infty}:=\lim_{t\to\infty}X_t$ exists almost surely
    and $1/X_{\infty}$ is a Bayes factor.\\
  \indent(2)
    Suppose $B$ is a precise Bayes factor.
    Then there is a test martingale $(X_t)$ such that $B=1/X_{\infty}$ a.s.
    Moreover, for any filtration $(\FFF_t)$
    such that $B$ is $\FFF_{\infty}$-measurable,
    there is a test martingale $(X_t,\FFF_t)$
    such that $B=1/X_{\infty}$ almost surely.
  \end{theorem}

\begin{pf}
  If $(X_t,\FFF_t)$ is a test supermartingale,
  the limit $X_{\infty}$ exists almost surely by Doob's convergence theorem
  (\cite{dellacheriemeyer1982}, VI.6),
  and the inequality\break $\int X_{\infty}\,\dd\Prob\le1$ holds by Fatou's lemma:
 \[
    \int X_{\infty}\,\dd\Prob
    =
    \int \liminf_{t\to\infty} X_t\,\dd\Prob
    \le
    \liminf_{t\to\infty} \int X_t\,\dd\Prob
    \le
    1.
 \]

  Now suppose that $B$ is a precise Bayes factor
  and $(\FFF_t)$ is a filtration
  (not necessarily satisfying the usual conditions)
  such that $B$ is $\FFF_{\infty}$-measurable;
  for concreteness, we consider the case of continuous time.
  Define a~test martingale $(X_t,\FFF^{\Prob}_{t+})$
  by setting $X_t:=\Expect(1/B\givn\FFF^{\Prob}_{t+})$;
  versions of conditional expectations can be chosen in such a way
  that $(X_t)$ is right-continuous:
  cf. \cite{dellacheriemeyer1982}, VI.4.
  Then $X_{\infty}=1/B$ almost surely by L\'{e}vy's zero--one law
  (\cite{levy1937}, pages 128--130; \cite{meyer1966}, VI.6, corollary).
  It remains to notice that $(X_t,\FFF_t)$ will also be a~test martingale.
  If $(\FFF_t)$ such that $B$ is $\FFF_{\infty}$-measurable
  is not given in advance, we can define it by, for example,
  \[
    \FFF_t
    :=
    \cases{
      \{\emptyset,\Omega\}, & if  $t<1$,\cr
      \sigma(B), & otherwise,
   }
 \]
  where $\sigma(B)$ is the $\sigma$-algebra generated by $B$.
\end{pf}

Formally, a \textit{stopping time} with respect to a filtration
$(\FFF_t)$ is a nonnegative random variable $\tau$ taking values in
$[0,\infty]$ such that, at each time $t$, the event
$\{\omega\st\tau(\omega)\le t\}$ belongs to $\FFF_t$. Let
$(X_t,\FFF_t)$ be a~test supermartingale. Doob's convergence theorem,
which was used in the proof of Theorem \ref{thm:bayes}, implies that we
can define its value $X_{\tau}$ at $\tau$ by the formula
$X_{\tau}(\omega):=X_{\tau(\omega)}(\omega)$ even when $\tau=\infty$
with positive probability. The \textit{stopped process}
$(X^{\tau}_t,\FFF_t):=(X_{t\wedge\tau},\FFF_t)$, where $a\wedge
b:=\min(a,b)$, will also be a~test supermartingale
(\cite{dellacheriemeyer1982}, VI.12). Since $X_{\tau}$ is the final
value of the stopped process, it follows from part 1 of
Theorem \ref{thm:bayes} that $1/X_{\tau}$ is a Bayes factor. (This also
follows directly from Doob's stopping theorem, \cite{meyer1966},
VI.13.)


\section{Supermartingales and $\lowercase{p}$-Values}
\label{sec:pvalues}

Now we will prove that the inverse of a supremum of a test
supermartingale is a $p$-test. This is true when the supremum is taken
over $[0,t]$ for some time point $t$ or over $[0,\tau]$ for any
stopping time $\tau$, but the strongest way of making the point is to
consider the supremum over all time points (i.e., for $\tau:=\infty$).

We will also show how to
construct a test martingale that has the inverse of a given $p$-test as
its supremum. Because the time scale for this martingale is artificial,
the value of the construction is more mathematical than directly
practical; it will help us prove Theorem \ref{thm:reduce} in
Section \ref{sec:reduce}. But it may be worthwhile to give an intuitive
explanation of the construction. This is easiest when the $p$-test has
discrete levels, because then we merely construct a sequence of bets.
Consider a $p$-test $p$ that is equal to $1$ with probability $1/2$, to
$1/2$ with probability $1/4$, to $1/4$ with probability $1/8$, etc.:
\[
  \Prob\{p = 2^{-n}\} = 2^{-n-1}
\]
for $n=0,1,\ldots.$ To see that a function on
$\Omega$ that takes these values with these probabilities is a $p$-test,
notice that when $2^{-n} \le \delta < 2^{-n+1}$,
\[
  \Prob \{ p \le \delta\} = \Prob\{p\le 2^{-n}\} = 2^{-n} \le \delta.
\]
Suppose that we learn first whether $p$ is $1$. Then, if it is not $1$,
we learn whether it is $1/2$. Then, if it is not $1/2$, whether it is
$1/4$, etc. To create the test martingale $X_0,X_1,\ldots,$ we start
with capital $X_0=1$ and bet it all against $p$ being $1$. If we lose,
$X_1=0$ and we stop. If we win, $X_1=2$, and we bet it all against $p$
being $1/2$, etc. Each time we have even chances of doubling our money
or losing it all. If $p=2^{-n}$, then our last bet will be against
$p=2^{-n}$, and the amount we will lose, $2^n$, will be $X^*_{\infty}$.
So $1/X^*_{\infty}=p$, as desired.

Here is our formal result:
\begin{theorem}\label{thm:pvalues}
(1)  If $(X_t,\FFF_t)$ is a test supermar\-tingale,
  $1/X^*_{\infty}$ is a $p$-test.\\
\indent(2)
  If $p$ is a precise $p$-test,
  there is a test mar\-tingale $(X_t)$
  such that $p=1/X_{\infty}^*$.
\end{theorem}

\begin{pf}
  The inequality
  $
    \Prob
    \{
      1/X_{\infty}^* \le \delta
    \}
    \le
    \delta
  $\break for test supermartingales
  follows from the maximal inequality~(\ref{eq:maximal}).

  In the opposite direction, let $p$ be a precise $p$-test.
  Set $\Pi:=1/p$;
  this function takes values in $[1,\infty]$.
  Define a right-continuous random process $(X_t)$, $t\in[0,\infty)$, by
  \[
    X_t(\omega)
    =
    \cases{
      1, & if  $t\in[0,1)$,\cr
      t, & if  $t\in[1,\Pi(\omega))$,\cr
      0, & otherwise.
    }
  \]
  Since $X_{\infty}^*=\Pi$,
  it suffices to check that $(X_t)$ is a~test martingale.
  The time interval where this process is nontrivial is $t\ge1$;
  notice that $X_1=1$ with probability one.

  Let $t\ge1$; we then have $X_t=t\III_{\{\Pi>t\}}$.
  Since $X_t$ takes values in the two-element set $\{0,t\}$,
  it is integrable.
  The $\sigma$-algebra generated by $X_t$ consists of 4 elements
  ($\emptyset$, $\Omega$, the set $\Pi^{-1}((t,\infty])$, and its complement),
  and the $\sigma$-algebra $\FFF_t$ generated by $X_s$, $s\le t$,
  consists of the sets $\Pi^{-1}(E)$
  where $E$ is either a~Borel subset of $[1,t]$
  or the union of $(t,\infty]$ and a~Borel subset of $[1,t]$.
  To check (\ref{eq:martingale}), where $1\le s<t$,
  it suffices to show that
  \[
    \int_{\Pi^{-1}(E)}
    X_t\,
    \dd\Prob
    =
    \int_{\Pi^{-1}(E)}
    X_s\,
    \dd\Prob,
  \]
  that is,
  \begin{equation}\label{eq:to-check}
    \int_{\Pi^{-1}(E)}
    t \III_{\{\Pi>t\}}\,
    \dd\Prob
    =
    \int_{\Pi^{-1}(E)}
    s \III_{\{\Pi>s\}}\,
    \dd\Prob,
  \end{equation}
  where $E$ is either a Borel subset of $[1,s]$
  or the union of $(s,\infty]$ and a Borel subset of $[1,s]$.
  If $E$ is a Borel subset of $[1,s]$,
  the equality (\ref{eq:to-check}) holds, as its two sides are zero.
  If $E$ is the union of $(s,\infty]$ and a Borel subset of $[1,s]$,
  (\ref{eq:to-check}) can be rewritten as
  \[
    \int_{\Pi^{-1}((s,\infty])}
    t \III_{\{\Pi>t\}}\,
    \dd\Prob
    =
    \int_{\Pi^{-1}((s,\infty])}
    s \III_{\{\Pi>s\}}\,
    \dd\Prob,
  \]
  that is,
  $
    t \Prob\{\Pi>t\}
    =
    s \Prob\{\Pi>s\}
  $,
  that is, $1=1$.
\end{pf}

\section{Calibrating $\lowercase{p}$-Values}
\label{sec:calibrate}

An increasing (not necessarily strictly increasing) function
$f\dvtx[0,1]\to[0,\infty]$ is called a \textit{calibrator} if $f(p)$ is a
Bayes factor for any $p$-test $p$. This notion was discussed in
\cite{vovk1993} and, less explicitly, in \cite{sellkeetal2001}. In this
section we will characterize the set of all increasing functions that
are calibrators; this result is a slightly more precise version of
Theorem 7 in \cite{vovk1993}.

We say that a calibrator $f$ \textit{dominates} a calibrator $g$ if
$f(x)\le g(x)$ for all $x\in[0,1]$. We say that $f$ \textit{strictly
dominates} $g$ if $f$ dominates $g$ and $f(x)<g(x)$ for some
$x\in[0,1]$. A calibrator is \textit{admissible} if it is not strictly
dominated by any other calibrator.

\begin{theorem}\label{thm:calibrate}
    (1)
    An increasing function $f\dvtx\break [0,1]\to[0,\infty]$ is a calibrator
    if and only if
    \begin{equation}\label{eq:calibrate-le}
      \int_0^1
      \frac{\dd x}{f(x)}
      \le
      1.
    \end{equation}\\
\indent(2)
    Any calibrator is dominated by an admissible calibrator.\\
\indent(3)
    A calibrator is admissible if and only if
    it is left-continuous and
    \begin{equation}\label{eq:calibrate-eq}
      \int_0^1
      \frac{\dd x}{f(x)}
      =
      1.
    \end{equation}
\end{theorem}

\begin{pf}
  Part 1 is proven in \cite{vovk1993} (Theorem 7),
  but we will give another argument, perhaps more intuitive.
  The condition ``only if'' is obvious:
  every calibrator must satisfy (\ref{eq:calibrate-le})
  in order to transform the ``exemplary'' $p$-test $p(\omega)=\omega$
  on the probability space $([0,1],\FFF,\Prob)$,
  where $\FFF$ is the Borel $\sigma$-algebra on $[0,1]$
  and $\Prob$ is the uniform probability measure on $\FFF$,
  into a Bayes factor.
  To check ``if,''
  suppose (\ref{eq:calibrate-le}) holds and take any $p$-test $p$.
  The expectation $\Expect(1/ f(p))$ depends on $p$
  only via the values $\Prob\{p\le c\}$, $c\in[0,1]$,
  and this dependence is monotonic:
  if a $p$-test $p_1$ is \textit{stochastically smaller} than another $p$-test $p_2$
  in the sense that $\Prob\{p_1\le c\} \ge \Prob\{p_2\le c\}$ for all $c$,
  then $\Expect(1/f(p_1))\ge\Expect(1/f(p_2))$.
  This can be seen, for example, from the well-known formula
  $\Expect(\xi)=\int_{0}^{\infty}\Prob\{\xi>c\}\,\dd c$,
  where $\xi$ is a nonnegative random variable:
  \begin{eqnarray*}
    \Expect\bigl(1/f(p_1)\bigr)
    &=&
    \int_{0}^{\infty}\Prob\{1/f(p_1)>c\}\,\dd c\\
    &\ge&
    \int_{0}^{\infty}\Prob\{1/f(p_2)>c\}\,\dd c
    =
    \Expect\bigl(1/f(p_2)\bigr).
  \end{eqnarray*}
  The condition (\ref{eq:calibrate-le}) means that the inequality\break $\Expect(1/ f(p))\le1$
  holds for our exemplary $p$-test $p$;
  since $p$ is stochastically smaller than any other $p$-test,
  this inequality holds for any $p$-test.

  Part 3 follows from part 1, and part 2 follows from parts 1 and 3.
\end{pf}

%

Equation (\ref{eq:calibrate-eq}) gives a recipe for producing
admissible calibrators $f$: take any left-continuous decreasing
function $g\dvtx [0,1]\to[0,\infty]$ such that $\int_0^1 g(x)\,\dd x = 1$ and
set $f(x):=1/g(x)$, $x\in[0,1]$. We see in this way, for example, that
%
\begin{equation}\label{eq:class-1}
  f(x)
  :=
  x^{1-\alpha}/\alpha
\end{equation}
is an admissible calibrator for every $\alpha\in(0,1)$; if we are
primarily interested in the behavior of $f(x)$ as $x\to0$, we should
take a small value of $\alpha$. This class of calibrators was found
independently in \cite{vovk1993} and~\cite{sellkeetal2001}.

The calibrators (\ref{eq:class-1}) shrink to 0 significantly slower
than $x$ as $x\to0$. But there are evidently calibrators that shrink as
fast as $x\ln^{1+\alpha}(1/x)$, or $x \ln(1/x)\cdot \ln^{1+\alpha}\ln(1/x)$,
etc., where $\alpha$ is a positive constant. For example,
\begin{equation}\label{eq:class-2}
  f(x)
  :=
  \cases{
    \alpha^{-1} (1+\alpha)^{-\alpha} x \ln^{1+\alpha}(1/x)\cr
     \hspace*{27.6pt}\mbox{if } x \le e^{-1-\alpha},\vspace*{2pt}\cr
    \infty, \quad\mbox{otherwise},
 }
\end{equation}
is an admissible calibrator for any $\alpha>0$.

%
%

\section{Calibrating the Running Suprema of Test Supermartingales}
\label{sec:reduce}

Let us call an increasing function $f\dvtx[1,\infty)\to[0,\infty)$ a
\textit{martingale calibrator} if it satisfies the following property:

\begin{quote}
  For any probability space $(\Omega,\FFF,\Prob)$
  and any test supermartingale $(X_t,\FFF_t)$ in this probability space
  there exists a test supermartingale $(Y_t,\FFF_t)$
  such that $Y_t\ge f(X^*_t)$ for all $t$ almost surely.
\end{quote}
There are at least 32 equivalent definitions of a martingale
calibrator: we can independently replace each of the two entries of
``supermartingale'' in the definition by ``martingale,'' we can
independently replace $(X_t,\FFF_t)$ by $(X_t)$ and $(Y_t,\FFF_t)$ by
$(Y_t)$, and we can optionally allow $t$ to take value $\infty$. The
equivalence will be demonstrated in the proof of
Theorem \ref{thm:reduce}.
Our convention is that $f(\infty):=\lim_{x\to\infty}f(x)$ (but remember
that $X^*_t=\infty$ only with probability zero, even for $t=\infty$).


As in the case of calibrators, we say that a martingale calibrator $f$
is \textit{admissible} if there is no other martingale calibrator $g$
such that $g(x)\ge f(x)$ for all $x\in[1,\infty)$ ($g$ \textit{dominates}
$f$) and $g(x)>f(x)$ for some $x\in[1,\infty)$.

\begin{theorem}\label{thm:reduce}
 (1) An increasing function $f\dvtx\break [1,\infty)\to[0,\infty)$ is a martingale calibrator
  if and only if
  \begin{equation}\label{eq:reduce-le}
    \int_0^1
    f(1/x)\,
    \dd x
    \le
    1.
  \end{equation}\\
\indent(2)
  Any martingale calibrator is dominated by an admissible martingale
  calibrator.\\
\indent(3)
  A martingale calibrator is admissible if and only if
  it is right-continuous and
  \begin{equation}\label{eq:reduce-eq}
    \int_0^1
    f(1/x)\,
    \dd x
    =
    1.
  \end{equation}
\end{theorem}

\begin{pf}
  We start from the statement ``if'' of\break part~1.
  Suppose an increasing function $f\dvtx[1,\infty)\to[0,\infty)$
  satisfies (\ref{eq:reduce-le})
  and $(X_t,\FFF_t)$ is a test supermartingale.
  By Theorem \ref{thm:calibrate},
  $g(x):=1/f(1/x)$, $x\in[0,1]$, is a calibrator,
  and by Theorem \ref{thm:pvalues}, $1/X^*_{\infty}$ is a $p$-test.
  Therefore,
  $g(1/X^*_{\infty})=1/f(X^*_{\infty})$
  is a Bayes factor,
  that is, $\Expect(f(X^*_{\infty}))\le1$.
  Similarly to the proof of Theorem \ref{thm:bayes},
  we set $Y_t := \Expect(f(X^*_{\infty})\givn\FFF_t)$,
  obtaining a nonnegative martingale $(Y_t,\FFF_t)$
  satisfying $Y_{\infty} = f(X^*_{\infty})$ a.s.
  We have $\Expect(Y_0)\le1$;
  the case $\Expect(Y_0)=0$ is trivial,
  and so we assume $\Expect(Y_0)>0$.
  Since
  \[
    Y_t
    =
    \Expect(f(X^*_{\infty})\givn\FFF_t)
    \ge
    \Expect(f(X^*_t)\givn\FFF_t)
    =
    f(X^*_t)\quad
    \mbox{a.s.}
  \]
  (the case $t=\infty$ was considered separately)
  and we can make $(Y_t,\FFF_t)$ a test martingale
  by dividing each $Y_t$ by $\Expect(Y_0)\in(0,1]$,
  the statement ``if'' in part 1 of the theorem is proven.
  Notice that our argument shows that $f$ is a martingale calibrator in any of the 32 senses;
  this uses the fact that $(Y_t)$ is a test (super)martingale
  whenever $(Y_t,\FFF_t)$ is a test (super)martingale.

  Let us now check that any martingale calibrator (in any of the senses)
  satisfies (\ref{eq:reduce-le}).
  By any of our definitions of a martingale calibrator,
  we have $\int f(X_t^*)\,\dd\Prob \le 1$
  for all test martingales $(X_t)$ and all $t<\infty$.
  It is easy to see that in Theorem \ref{thm:pvalues}, part 2,
  we can replace $X^*_{\infty}$ with, say, $X^*_{\pi/2}$,
  by replacing the test martingale $(X_t)$ whose existence it asserts with
 \[
    X'_t
    :=
    \cases{
      X_{\tan t}, & if $t<\pi/2$,\cr
      X_{\infty}, & otherwise.
   }
  \]
  Applying this modification of Theorem \ref{thm:pvalues}, part 2,
  to the precise $p$-test $p(\omega):=\omega$
  on $[0,1]$ equipped with the uniform probability measure,
  we obtain
  \[
    1
    \ge
    \int f(X_{\pi/2}^*)\,\dd\Prob
    =
    \int f(1/p)\,\dd\Prob
    =
    \int_0^1 f(1/x)\,\dd x.
  \]
  This completes the proof of part 1.

  Part 3 is now obvious, and part 2 follows from parts~1 and 3.
\end{pf}

As in the case of calibrators, we have a recipe for producing
admissible martingale calibrators $f$ provided by (\ref{eq:reduce-eq}):
take any left-continuous decreasing function $g\dvtx[0,1]\to[0,\infty)$
satisfying $\int_0^1 g(x)\,\dd x=1$ and set $f(y):=g(1/y)$,
$y\in[1,\infty)$. In this way we obtain the class of admissible
martingale calibrators
%
\begin{equation}\label{eq:martingale-class-1}
  f(y)
  :=
  \alpha y^{1-\alpha},\quad
  \alpha\in(0,1),
\end{equation}
analogous to (\ref{eq:class-1}) and the class
\[
  f(y)
  :=
  \cases{
    \alpha (1+\alpha)^{\alpha} \dfrac{y}{\ln^{1+\alpha}y},
      & if $y \ge e^{1+\alpha}$,\cr
    0,  & otherwise,
 }
  \quad\alpha>0,
\]
analogous to (\ref{eq:class-2}).

In the case of discrete time, Theorem \ref{thm:reduce} has been greatly
generalized by Dawid et al. (\cite{dawidetal2011}, Theorem~1). The
generalization, which required new proof techniques, makes it possible
to apply the result in new fields, such as mathematical finance
(\cite{dawidetal2011}, Section~4). 

In this article we have considered only tests of simple statistical
hypotheses. We can use similar ideas for testing composite hypotheses,
that is, sets of probability measures. One possibility is to measure
the evidence against the composite hypothesis by the current value of a
random process that is a~test supermartingale under all probability
measures in the composite hypothesis; we will call such processes
\textit{simultaneous test supermartingales}. For example,\break there are
nontrivial processes that are test supermartingales under all
exchangeable probability measures simultaneously (\cite{vovketal2005},
Section 7.1). Will martingale calibrators achieve their goal for
simultaneous test supermartingales? The method of proof of
Theorem \ref{thm:reduce} does not work in this situation: in general,
it will produce a different test supermartingale for each probability
measure. The advantage of the method used in \cite{dawidetal2011} is
that it will produce one process, thus demonstrating that for each
martingale calibrator $f$ and each simultaneous test supermartingale
$X_t$ there exists a simultaneous test supermartingale $Y_t$ such that
$Y_t\ge f(X_t^*)$ for all $t$ (the method of \cite{dawidetal2011} works
pathwise and makes the qualification ``almost surely'' superfluous).



\section{Examples}
\label{sec:example}

Although our results are very general, we can illustrate them using the
simple problem of testing whether a coin is fair. Formally, suppose we
observe a sequence of independent identically distributed binary random
variables $x_1,x_2,\ldots,$ each taking values in the set $\{0,1\}$;
the probability $\theta\in[0,1]$ of $x_1=1$ is unknown. Let
$P_{\theta}$ be the probability distribution of $x_1,x_2,\ldots$; it is
a probability measure on $\{0,1\}^{\infty}$. In most of this section,
our null hypothesis is that $\theta=1/2$.

\begin{figure*}[b]

\includegraphics{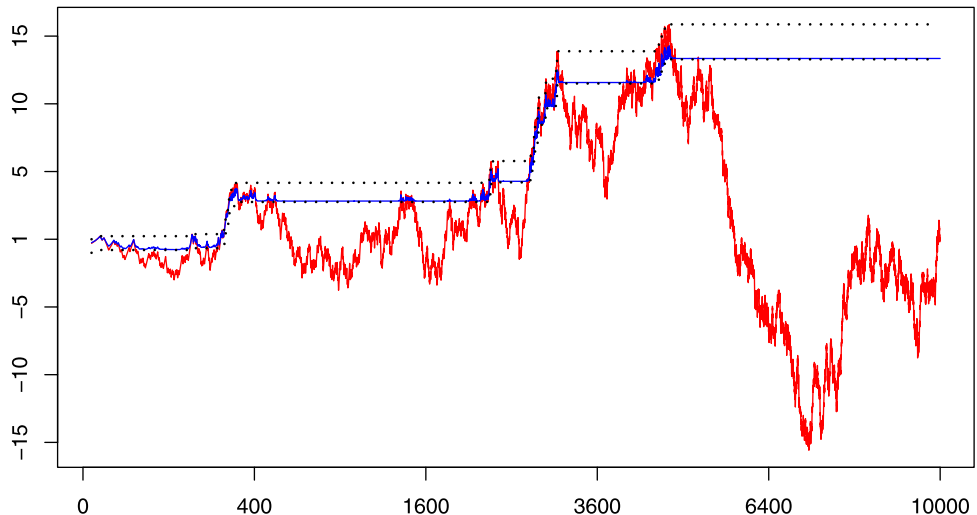}

    \caption{The red line is a realization over $10{,}000$ trials
      of the likelihood ratio for testing $\theta = 1/2$ against $\theta=3/4$.
      The horizontal axis gives the number of observations so far.
      The vertical axis is logarithmic and is labeled by powers of $10$.
      The likelihood ratio varies wildly, up to~$10^{15}$ and down to $10^{-15}$.
      Were the sequence continued indefinitely,
      it would be unbounded in both directions.}\label{fig:ss}
  \end{figure*}

We consider both Bayesian testing of $\theta=1/2$, whe\-re the output is
a posterior distribution, and non-Bayesian testing, where the output is
a $p$-value. We call the approach that produces $p$-values the
\textit{samp\-ling-theory approach} rather than the frequentist approach,
because it does not require us to interpret all probabilities as
frequencies; instead, we can merely interpret the $p$-values using
Cournot's principle (\cite{shafer2006}, Section 2). We have borrowed
the term ``sampling-theory'' from D. R. Cox and A. P. Dempster
\cite{cox2006,dempster1969}, without necessarily using it in exactly
the same way as either of them do.

We consider two tests of $\theta=1/2$, corresponding to two different
alternative hypotheses:
\begin{longlist}[(2)]
\item[(1)]
  First, we test $\theta=1/2$ against $\theta=3/4$.
  This is unrealistic on its face;
  it is hard to imagine accepting a~model that contains only these two simple hypotheses.
  But some of what we learn from this test will carry over
  to sensible and widely used tests of a simple against a composite hypothesis.
\item[(2)]
  Second, we test $\theta=1/2$ against the composite hypothesis $\theta\ne1/2$.
  In the spirit of Bayesian statistics and following Laplace
  (\cite{laplace1774};
  see also \cite{todhunter1865}, Section~870, and \cite{stigler1986}),
  we represent this composite hypothesis by the uniform distribution on $[0,1]$,
  the range of possible values for $\theta$.
  (In general, the composite hypotheses of this section
  will be composite\vadjust{\goodbreak} only in the sense of Bayesian statistics;
  from the point of view of the sampling-theory approach,
  these are still simple hypotheses.)
\end{longlist}
For each test, we give an example of calibration of the running
supremum of the likelihood ratio. In the case of the composite
alternative hypothesis, we also discuss the implications of using the
inverse of the running supremum of the likelihood ratio as
a $p$-value.

To round out the picture, we also discuss Bayesian testing of the
composite hypothesis $\theta\le1/2$ against the composite hypothesis
$\theta>1/2$, representing the former by the uniform distribution on
$[0,1/2]$ and the latter by the uniform distribution on $(1/2,1]$.
Then, to conclude, we discuss the relevance of the calibration of
running suprema to Bayesian philosophy.

Because the idea of tracking the supremum of a~martingale is related to
the idea of waiting until it reaches a~high value, our discussion is
related to a long-standing debate about ``sampling to reach a foregone
conclusion,'' that is, continuing to sample in search of evidence
against a hypothesis and stopping only when some conventional $p$-value
finally dips below a conventional level such as $5\%$. This debate goes
back at least to the work of Francis Anscombe in 1954
\cite{anscombe1954}. In 1961, Peter Armitage described situations where
even a~Bayesian can sample to a foregone conclusion
(\cite{armitage1961}; \cite{bernardosmith2000}, Section 5.1.4). Yet in
1963 \cite{edwardsetal1963}, Ward Edwards and his co-authors insisted
that this is not a problem: ``The likelihood principle emphasized in
Bayesian statistics implies, among other things, that the rules
governing when data collection stops are irrelevant to data
interpretation. It is entirely appropriate to collect data until a
point has been proven or disproven, or until the data collector runs
out of time, money, or patience.'' For further information on this
debate, see \cite{wagenmakers2007}. We will not attempt to analyze it
thoroughly, but our examples may be considered a contribution to it.

\subsection{Testing $\theta=1/2$ Against a Simple Alternative}
\label{subsec:ss}

To test our null hypothesis $\theta=1/2$ against the alternative
hypothesis $\theta = 3/4$, we use the likelihood ratio
\begin{eqnarray}\label{eq:simple-process}
  X_t
  &:=&
  \frac{P_{3/4}(x_1,\ldots,x_t)}{P_{1/2}(x_1,\ldots,x_t)}\nonumber\\ [-8pt]\\ [-8pt]
 & =&
  \frac{(3/4)^{k_t}(1/4)^{t-k_t}}{(1/2)^t}
  =
  \frac{3^{k_t}}{2^t},\nonumber
\end{eqnarray}
where $k_t$ is the number of 1s in $x_1,\ldots,x_t$ [and
$P_{\theta}(x_1,\break\ldots,x_t)$ is the probability under $P_\theta$ that
the first $t$ observations are $x_1,\ldots,x_t$; such informal notation
was already used in Section \ref{sec:history}]. The sequence of
successive values of this likelihood ratio is a~test martingale
$(X_t)$.

According to (\ref{eq:martingale-class-1}), the function
%
\begin{equation}\label{eq:example5}
  f(y)
  :=
  0.1 y^{0.9}
\end{equation}
is a martingale calibrator. So there exists a test martingale $(Y_t)$
such that
%
\begin{equation}\label{eq:calibrating-ratio}
  Y_t \ge \max_{n=1,\ldots,t}  0.1 X_n^{0.9}.
\end{equation}


Figure \ref{fig:ss} shows an example in which the martingale
calibrator (\ref{eq:example5}) preserves\vadjust{\goodbreak} a reasonable amount of the
evidence against $\theta=1/2$. To construct this figure, we generated a
sequence $x_1,\ldots,x_{10,000}$ of $0$s and $1$s, choosing each $x_t$
independently with the probability $\theta$ for $x_t = 1$ always equal
to $\ln 2/ \ln 3 \approx 0.63$. Then we formed the lines in the figure
as follows:
\begin{itemize}
\item
  The red line is traced by the sequence of numbers
  $
    X_t = 3^{k_t}/2^t
  $.
  If our null hypothesis $\theta=1/2$ were true,
  these numbers would be a realization of a~test martingale,
  but this hypothesis is false (as is our alternative hypothesis $\theta = 3/4$).
\item
  The upper dotted line is the running supremum of the $X_t$:
  \begin{eqnarray*}
    X^*_t  &=&  \max_{n=1,\ldots,t}  \frac{3^{k_n}}{2^n}\\
           &=& \mbox{(best evidence so far against $\theta=1/2$)}_t.
  \end{eqnarray*}
\item
  The lower dotted line, which we will call $F_t$,
  shrinks this best evidence using our martingale calibrator:
  $
    F_t = 0.1(X^*_t )^{0.9}
  $.
\item
  The blue line, which we will call $Y_t$, is a test martingale
  under the null hypothesis that satisfies (\ref{eq:calibrating-ratio}):
  $
    Y_t \ge F_t
  $.
\end{itemize}
According to the proof of Theorem \ref{thm:reduce},
$\Expect(0.1(X^*_{\infty})^{0.9}\givn\break
\FFF_t)/\Expect(0.1(X^*_{\infty})^{0.9})$, where the expected values
are with respect to $P_{1/2}$, is a test martingale that satisfies
(\ref{eq:calibrating-ratio}). Because these expected values may be
difficult to compute, we have used in its stead in the role of $Y_t$
a~more easily computed test martingale that is shown in
\cite{dawidetal2011} to satisfy (\ref{eq:calibrating-ratio}).

Here are the final values of the processes shown in
Figure \ref{fig:ss}:
\begin{eqnarray*}
  X_{10,000} &=& 2.2, \quad \hspace*{30.2pt}X^*_{10,000} = 7.3 \times 10^{15}, \\
  F_{10,000} &=& 1.9 \times 10^{13},\quad  Y_{10,000} = 2.2 \times 10^{13}.
\end{eqnarray*}
The test martingale $Y_t$ legitimately and correctly rejects the null
hypothesis at time $10{,}000$ on the basis of $X_t$'s high earlier
values, even though the Bayes factor $X_{10{,}000}$ is not high. The
Bayes factor $Y_{10,000}$ gives overwhelming evidence against the null
hypothesis, even though it is more than two orders of magnitude smaller
than\vspace*{1pt} $X^*_{10,000}$.

\begin{figure*}[b]

\includegraphics{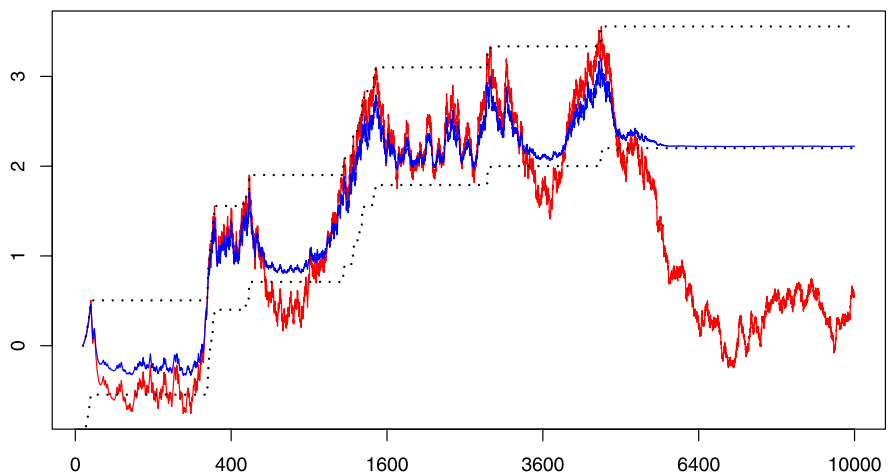}

    \caption{A realization over $10{,}000$ trials of the likelihood ratio
      for testing $\theta = 1/2$ against
      the probability distribution $Q$ obtained by averaging $P_{\theta}$
      with respect to the uniform distribution for $\theta$.
      The vertical axis is again logarithmic.
      As in Figure~\textup{\protect\ref{fig:ss}},
      the oscillations would be unbounded
      if trials continued indefinitely.}
      \label{fig:cs}
  \end{figure*}

As the reader will have noticed, the test martingale $X_t$'s
overwhelming values against $\theta=1/2$ in Figure \ref{fig:ss} are
followed, around $t=7{,}000$, by overwhelming values (order of magnitude
$10^{-15}$) against $\theta=3/4$. Had we been testing $\theta=3/4$
against $\theta=1/2$, we would have found that it can also be rejected
very strongly even after calibration. The fact that $(X_t)$  and
$(1/X_t)$ both have times when they are very large is not accidental
when we sample from $P_{\ln 2/ \ln 3}$. Under this measure, the
conditional expected value of the increment $\ln X_t-\ln X_{t-1}$,
given the first $t-1$ observations, is
\[
  \frac{\ln2}{\ln3}
  \ln\frac32
  +
  \biggl(1-\frac{\ln2}{\ln3}\biggr)
  \ln\frac12
  =
  0.
\]
So $\ln X_t$ is a martingale under $P_{\ln2/\ln3}$. The conditional
variance of its increment is
\[
  \frac{\ln2}{\ln3}
  \biggl(\ln\frac32\biggr)^2
  +\biggl(1 - \frac{\ln2}{\ln3}\biggr)
  \biggl(\ln\frac12\biggr)^2
  =
  \ln2\ln\frac32.
\]
By the law of the iterated logarithm,
\[
  \limsup_{t\to\infty}
  \frac{\ln X_t}
  {\sqrt{2\ln2\ln(3/2)  t\ln\ln t}}
  =
  1
\]
 and
 \[
  \liminf_{t\to\infty}
  \frac{\ln X_t}
  {\sqrt{2\ln2\ln(3/2)  t\ln\ln t}}
  =
  -1
\]
almost surely. This means that as $t$ tends to $\infty$, $\ln X_t$
oscillates between approximately $\pm0.75\sqrt{t\ln\ln t}$; in
particular,
%
\begin{equation}\label{eq:osc}
  \limsup_{t\to\infty}X_t=\infty
  \quad\mbox{and}\quad
  \liminf_{t\to\infty}X_t=0
\end{equation}
almost surely. This guarantees that we will eventually obtain
overwhelming evidence against whichever of the hypotheses $\theta=1/2$
and $\theta=3/4$ that we want to reject.  This may be called sampling
to a foregone conclusion, but the foregone conclusion will be correct,
since both $\theta=1/2$ and $\theta=3/4$ are wrong.

In order to obtain (\ref{eq:osc}), we chose $x_1,\ldots,x_{10,000}$
from a probability distribution, $P_{\ln 2/\ln3}$, that lies midway
between $P_{1/2}$ and $P_{3/4}$ in the sense that it tends to produce
sequences that are as atypical with respect to the one measure as to
the other. Had we chosen a~sequence $x_1,\ldots,x_{10,000}$ less
atypical with respect to $P_{3/4}$ than with respect to $P_{1/2}$, then
we might have been able to sample to the foregone conclusion of
rejecting $\theta=1/2$, but not to the foregone conclusion of rejecting
$\theta=3/4$.

\subsection{Testing $\theta=1/2$ Against a Composite Alternative}
\label{subsec:sc}

Retaining $\theta=1/2$ as our null hypothesis, we now take as our
alternative hypothesis the probability distribution $Q$ obtained by
averaging $P_{\theta}$ with respect to the uniform distribution for
$\theta$.

\begin{figure*}

\includegraphics{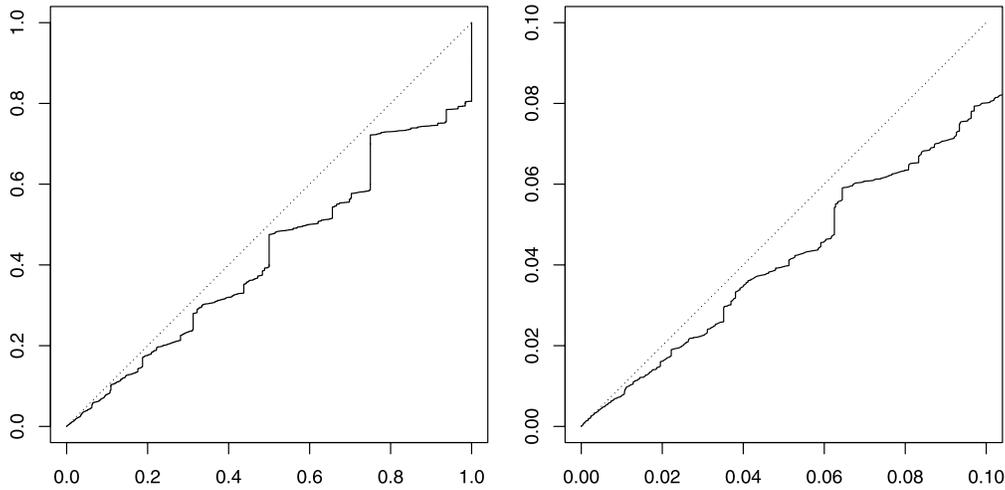}

    \caption{On the left we graph $\Prob\{p\le\delta\}$ as a function of $\delta$,
      where $p$ is the function defined in (\protect\ref{eq:$p$-test}).
      On the right, we magnify the lower left corner of this graph.}\label{fig:p}
  \end{figure*}

After we observe $x_1,\ldots,x_t$, the likelihood ratio for testing
$P_{1/2}$ against $Q$ is
%
\begin{eqnarray}\label{eq:process}
  \hspace*{20pt}X_t
  :&=&
  \frac{Q(x_1,\ldots,x_t)}{P_{1/2}(x_1,\ldots,x_t)}\nonumber\\ [-8pt]\\ [-8pt]
 \phantom{:}&=&
  \frac{\int_0^1\theta^{k_t}(1-\theta)^{t-k_t}\,\dd\theta}{(1/2)^t}
  =
  \frac{k_t!(t-k_t)!2^t}{(t+1)!}.\nonumber
\end{eqnarray}
Figure \ref{fig:cs} shows an example of this process and of the
application of the same martingale calibrator, (\ref{eq:example5}),
that we used in Figure \ref{fig:ss}. In this case, we generate the 0s
and 1s in the sequence $x_1,\ldots,x_{10,000}$ independently but with a
probability for $x_t=1$ that slowly converges to $1/2$:
$\frac12+\frac14\sqrt{\ln t/t}$. As we show in
Appendix \ref{app:calculations}, (\ref{eq:osc}) again holds almost
surely; if you wait long enough, you will have enough evidence to
reject legitimately whichever of the two false hypotheses
(independently and identically distributed with $\theta=1/2$, or
independently and identically distributed with $\theta\ne1/2$) you
want.

Here are the final values of the processes shown in
Figure \ref{fig:cs}:
\begin{eqnarray*}
  X_{10,000}, &=& 3.5,  \hspace*{11.2pt}X^*_{10,000} = 3599, \\
    F_{10,000}, &=& 159,\quad Y_{10,000} = 166.
\end{eqnarray*}
In this case, the evidence against $\theta=1/2$ is very substantial but
not overwhelming.

\vspace*{3pt}
\subsection{$p$-Values for Testing $\theta=1/2$}
\vspace*{3pt}

By Theorem \ref{thm:pvalues}, $1/X^*_{\infty}$ is a $p$-test whenever
$(X_t)$ is a test martingale.  Applying this to the test
martingale (\ref{eq:process}) for testing $P_{1/2}$ against $Q$, we see
that
%
\begin{eqnarray}\label{eq:$p$-test}
  p(x_1,x_2,\ldots)
  :&=&
  \frac{1}{\sup_{1\le t <
  \infty}(k_t!(t-k_t)!2^t/(t+1)!)}\nonumber\\[-7pt]\\ [-7pt]
  \phantom{:}&=&
  \inf_{1\le t < \infty}\frac{(t+1)!}{k_t!(t-k_t)!2^t}\nonumber
\end{eqnarray}
is a $p$-test for testing $\theta=1/2$ against $\theta\ne 1/2$. Figure~\ref{fig:p} shows that it is only moderately
conservative.\looseness=1

Any function of the observations that is bounded below by a $p$-test is
also a $p$-test. So for any rule $N$ for selecting a positive integer
$N(x_1,x_2,\ldots)$ based on knowledge of some or all of the
observations $x_1,\break x_2,\ldots,$ the function
%
\begin{equation}\label{eq:rN}
   r_N(x_1,x_2,\ldots) := \frac{(N+1)!}{k_N!(N-k_N)!2^N}
\end{equation}
is a $p$-test. It does not matter whether $N$ qualifies as a stopping
rule [i.e., whether $x_1,\ldots,x_n$ always determine whether
$N(x_1,x_2,\ldots)\le n$].

For each positive integer $n$, let
\begin{equation}\label{eq:pn}
     p_n := \frac{(n+1)!}{k_n!(n-k_n)!2^n}.
\end{equation}
We can paraphrase the preceding paragraph by saying that
$p_n$ is a $p$-value (i.e., the value of a $p$-test) no matter what rule is
used to select $n$. In parti\-cular, it is a~$p$-value even if it was
selected because it was the smallest number in the sequence
$p_1,p_2,\ldots,p_n,\ldots,p_t$, where~$t$ is an integer much larger
than $n$.

We must nevertheless be cautious if we do not know the rule $N$---if
the experimenter who does the sampling reports to us $p_n$ and perhaps
some other information but not the rule $N$. We can consider the
reported value of $p_n$ a legitimate $p$-value whenever we know that the
experimenter would have told us $p_n$ for some $n$, even if we do not
know what rule $N$ he followed to choose~$n$ and even if he did not
follow any clear rule. But we should not think of $p_n$ as a $p$-value if
it is possible that the experimenter would not have reported anything
at all had he not found an $n$ with a $p_n$ to his liking. We are
performing a $p$-test only if we learn the result no matter what it
is.\looseness=-1

\begin{figure*}

\includegraphics{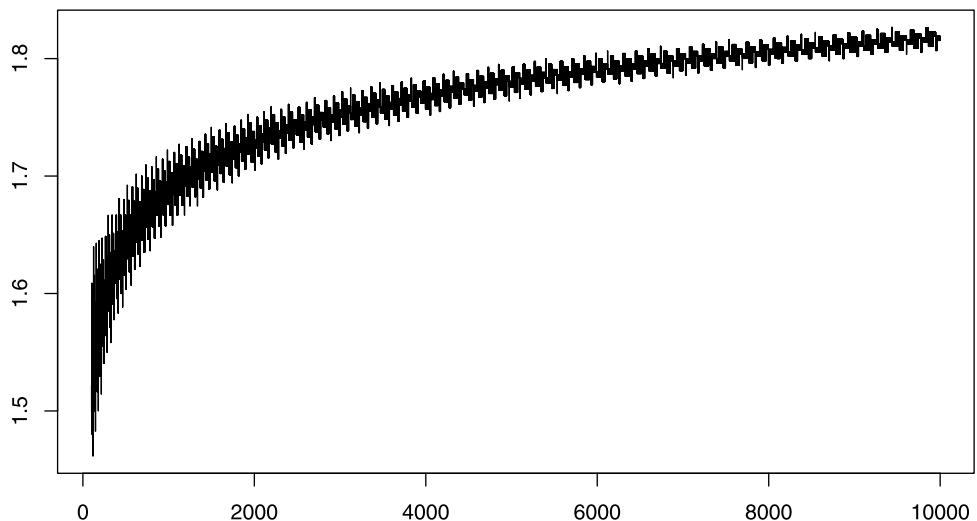}

    \caption{The ratio (\protect\ref{eq:width-ratio}) as $n$ ranges from $100$ to $10{,}000$.
      This is the factor by which not knowing $n$ in advance widens
      the $99\%$ prediction interval for $k_n$.
      Asymptotically, the ratio tends to infinity with $n$
      as $c\sqrt{\ln n}$ for some positive constant $c$.}\label{fig:comp}
  \end{figure*}

Continuing to sample in search of evidence against $\theta=1/2$ and
stopping only when the $p$-value finally reaches $5\%$ can be considered
legitimate if instead of using conventional $p$-tests for fixed sample
sizes we use the $p$-test (\ref{eq:rN}) with $N$ defined by
\[
  N(x_1,x_2,\ldots)
  :=
  \inf  \biggl\{ n \Bigm|  \frac{(n+1)!}{k_n!(n-k_n)!2^n} \le 0.05\biggr\}.
\]
But we must bear in mind that $N(x_1,x_2,\ldots)$ may take the value
$\infty$. If the experimenter stops only when the $p$-value dips down to
the $5\%$ level, he has a chance of at least $95\%$, under the null
hypothesis, of never stopping. So it will be legitimate to interpret a
reported $p_n$ of $0.05$ or less as a $p$-value (the observed value of a
$p$-test)
only if we were somehow also guaranteed to hear about the failure to
stop.

\subsection{Comparison with a Standard $p$-Test}

If the number $n$ of observations is known in advance, a standard
sampling-theory procedure for tes\-ting the hypothesis $\theta=1/2$ is to
reject it if $
  |k_n-n/2|\ge c_{n,\delta},
$ where $c_{n,\delta}$ is chosen so that
$P_{1/2}\{|k_n-n/2|\ge c_{n,\delta}\}$ is equal (or less
than but as close as possible) to a chosen significance level $\delta$.
To see how this compares with the $p$-value $p_n$ given by (\ref{eq:pn}),
let us compare the conditions for nonrejection:
\begin{itemize}
\item
  If we use the standard procedure,
  the condition for not rejecting $\theta=1/2$ at level $\delta$ is
  %
  \begin{equation}\label{eq:notrejs}
    |k_n-n/2|<c_{n,\delta}.
  \end{equation}
\item
  If we use the $p$-value $p_n$, the condition for not rejecting $\theta=1/2$
  at level $\delta$ is $p_n > \delta$, or
  %
  \begin{equation}\label{eq:notrejm}
    \frac{(n+1)!}{k_n!(n-k_n)!2^n} > \delta.
  \end{equation}
\end{itemize}
In both cases, $k_n$ satisfies the condition with probability at least
$1-\delta$ under the null hypothesis, and, hence, the condition defines
a level $1-\delta$ prediction interval for~$k_n$. Because
condition~(\ref{eq:notrejs}) requires the value of $n$ to be known in
advance and condition~(\ref{eq:notrejm}) does not, we can expect the
prediction interval defined by (\ref{eq:notrejm}) to be wider than the
one determined by (\ref{eq:notrejm}). How much wider?

Figure \ref{fig:comp} answers this question for the case where $\delta
= 0.01$ and $100 \le n \le 10{,}000$. It shows, for each value of $n$ in
this range, the ratio
%
\begin{eqnarray}\label{eq:width-ratio}
  \frac{\mbox{width of the $99\%$ prediction interval given by (\ref{eq:notrejm})}}
       {\mbox{width of the $99\%$ prediction interval given by
       (\ref{eq:notrejs})}},\hspace*{8pt}\nonumber\\
\end{eqnarray}
that is, the factor by which not knowing $n$ in advance widens the
prediction interval. The factor is less than $2$ over the whole range
but increases steadi\-ly with $n$.

As $n$ increases further, the factor by which the standard interval is
multiplied increases without limit, but very slowly. To verify this, we
first rewrite (\ref{eq:notrejm})~as
\begin{equation}\label{eq:variable}
  \hspace*{15pt}|k_n-n/2|
  <
  (1+\alpha_n)
  \sqrt{n}
  \sqrt{\frac12\ln\frac{1}{\delta}+\frac14\ln n},
\end{equation}
where $\alpha_n$ is a sequence such that $\alpha_n\to0$ as
$n\to\infty$. [For some $\alpha_n$ of order $o(1)$ the
inequality (\ref{eq:variable}) is stronger than $p_n>\delta$, whereas
for others it is weaker; see Appendix \ref{app:calculations} for
details of calculations.] Then, using the Berry--Esseen theorem and
letting $z_{\epsilon}$ stand for the upper $\epsilon$-quantile of the
standard Gaussian distribution, we rewrite (\ref{eq:notrejs}) as
%
\begin{equation}\label{eq:fixed}
  \biggl|k_n-\frac{n}{2}\biggr|
  <
  \frac12
  z_{\delta/2+\alpha_n}
  \sqrt{n},
\end{equation}
where $\alpha_n$ is a sequence such that
$|\alpha_n|\le\break (2\pi)^{-1/2}n^{-1/2}$ for all $n$. (See
\cite{hippmattner2007}.)
As $\delta\to0$,
\[ 
  z_{\delta/2}
  \sim
  \sqrt{2\ln\frac{2}{\delta}}
  \sim
  \sqrt{2\ln\frac{1}{\delta}}.
\]
%
%
So the main asymptotic difference between
(\ref{eq:variable}) and (\ref{eq:fixed}) is the presence of the term
$\frac14\ln n$ in (\ref{eq:variable}).

The ratio (\ref{eq:width-ratio}) tends to infinity with $n$ as
$c\sqrt{\ln n}$ for a positive constant $c$ (namely, for
$c=1/z_{\delta/2}$, where $\delta=0.01$ is the chosen significance
level). However, the expression on the right-hand side of
(\ref{eq:variable}) results from using the uniform probability measure
on $\theta$ to average the probability measures $P_{\theta}$. Averaging
with respect to a different probability measure would give something
different, but it is clear from the law of the iterated logarithm that
the best we can get
is a prediction interval whose ratio with the standard interval will
grow like $\sqrt{\ln \ln n}$ instead of $\sqrt{\ln n}$. In fact, the
method we just used to obtain (\ref{eq:variable}) was used by Ville,
with a more carefully chosen probability measure on $\theta$, to prove
the upper half of the law of the iterated logarithm (\cite{ville1939},
Section V.3), and Ville's argument was rediscovered
and simplified using the algorithmic theory of randomness in
\cite{vovk1987full}, Theorem 1.


\begin{figure*}

\includegraphics{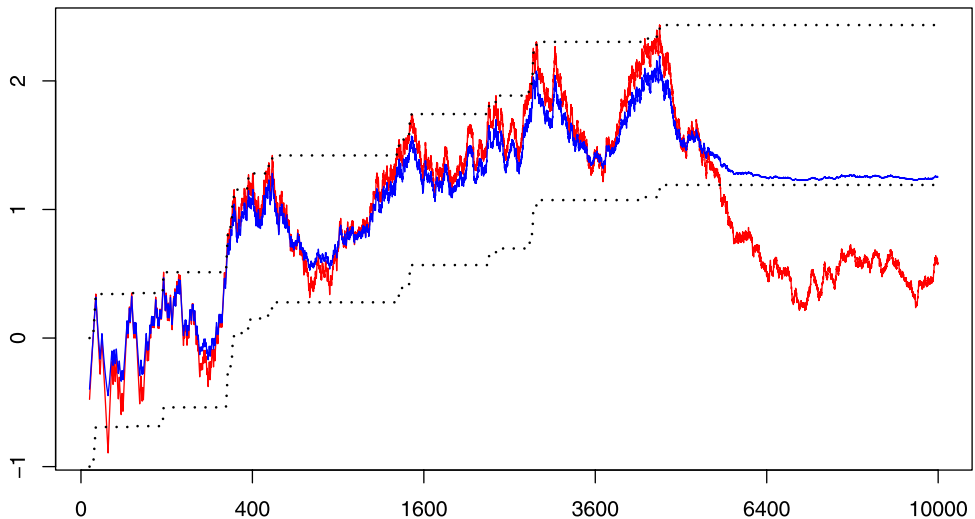}

    \caption{A realization over $10{,}000$ trials of the likelihood ratio
      for testing the probability distribution obtained by averaging $P_{\theta}$
      with respect to the uniform probability measure on $[0,1/2]$
      against the probability distribution obtained by averaging $P_{\theta}$
      with respect to the uniform probability measure on $(1/2,1]$.
      As in the previous figures, the vertical axis is logarithmic,
      and the red line would be unbounded in both directions
      if observations continued indefinitely.}\label{fig:cc}
      \vspace*{6pt}
  \end{figure*}

\subsection{Testing a Composite Hypothesis Against a~Composite Hypothesis}

When Peter Armitage pointed out that even Baye\-sians can sample to a
foregone conclusion, he used as an example the Gaussian model with
known variance and unknown mean \cite{armitage1961}. 
%
%
We can adapt Armita\-ge's idea to coin tossing by comparing
two composite hypotheses: the null hypothesis $\theta\le1/2$,
represented by the uniform probability measure on $[0,1/2]$, and the
alternative hypothesis $\theta>1/2$, represented by the uniform
probability measure on $(1/2,1]$. (These hypotheses are natural in the
context of paired comparison: see, e.g., \cite{lehmann2006}, Section
3.1.) The test martingale is
%
\begin{eqnarray}\label{eq:cc-1}
  X_t
 & =&
  \frac
  {2\int_{1/2}^1\theta^{k_t}(1-\theta)^{t-k_t}\,\dd\theta}
  {2\int_{0}^{1/2}\theta^{k_t}(1-\theta)^{t-k_t}\,\dd\theta}\nonumber\\
  [-8pt]\\ [-8pt]
  &=&
  \frac{\Prob\{B_{t+1}\le k_t\}}{\Prob\{B_{t+1}\ge k_t+1\}},\nonumber
\end{eqnarray}
where $B_n$ is the binomial random variable with parameters $n$ and
$1/2$; see Appendix \ref{app:calculations} for details. If the sequence
$x_1,x_2,\ldots$ turns out to be typical of $\theta=1/2$, then by the
law of the iterated logarithm, $(k_t-t/2)/\sqrt{t}$ will almost surely
have $\infty$ as its upper limit and $-\infty$ as its lower limit;
therefore, (\ref{eq:osc}) will hold again. This confirms Armitage's
intuition that arbitrarily strong evidence on both sides will emerge if
we wait long enough, but the oscillation depends on increasingly
extreme reversals of a random walk, and the lifetime of the universe
may not be long enough for us to see any of them
[$\sqrt{\ln\ln(5\times10^{23})}<2$]. 

Figure \ref{fig:cc} depicts one example, for which the final values are
\begin{eqnarray*}
  X_{10,000} &=& 3.7,\hspace*{13.5pt}  X^*_{10,000} = 272, \\
  F_{10,000} &=& 15.5, \quad Y_{10,000} = 17.9.
\end{eqnarray*}
In this realization, the first $10{,}000$ observations provi\-de modest
evidence against $\theta\le 1/2$ and none against $\theta>1/2$. Figures
\ref{fig:ss} and \ref{fig:cs} are reasonably typical for their setups,
but in this setup it is unusual for the first $10{,}000$ observations to
show even as much evidence against one of the hypotheses as we see in
Figure \ref{fig:cc}. 


\subsection{A Puzzle for Bayesians}

From a Bayesian point of view, it may seem puzzling that we should want
to shrink a likelihood ratio in order to avoid exaggerating the
evidence against a null hypothesis. Observations affect Bayesian
posterior odds only through the likelihood ratio, and we know that the
likelihood ratio is not affected by the sampling plan. So why should we
adjust it to take the sampling plan into account?

Suppose we assign equal prior probabilities of $1/2$ each to the two
hypotheses $\theta=1/2$ and $\theta=3/4$ in our first coin-tossing
example. Then if we stop at time $t$, the likelihood ratio $X_t$ given
by (\ref{eq:simple-process}) is identical with the posterior odds in
favor of $\theta=3/4$. If we write $\mathbf{post}_t$ for the posterior
probability measure at time $t$, then
\[
  X_t
  =
  \frac{\mathbf{post}_t\{\theta=3/4\}}{\mathbf{post}_t\{\theta=1/2\}}
  =
  \frac{1-\mathbf{post}_t\{\theta=1/2\}}{\mathbf{post}_t\{\theta =1/2\}}
\]
and
\begin{equation}\label{eq:post-model}
  \mathbf{post}_t\{\theta=1/2\} = \frac{1}{X_t+1}.
\end{equation}
This is our posterior probability given the evidence $x_1,\ldots,x_t$
no matter why we decided to stop at time~$t$. If we ``calibrate'' $X_t$
and plug the calibrated value instead of the actual value
into (\ref{eq:post-model}), we will get the posterior probability
wrong.

It may help us escape from our puzzlement to acknowledge that if the
model is wrong, then the observations may oscillate between providing
overwhelming evidence against $\theta=1/2$ and providing overwhelming
evidence against $\theta=3/4$, as in Figure \ref{fig:ss}. Only if we
insist on retaining the model in spite of this very anomalous
phenomenon will (\ref{eq:post-model}) continue to be our posterior
probability for $\theta=1/2$ at time $t$, and it is this stubbornness
that opens the door to sampling to whichever foregone conclusion we
want, $\theta=1/2$ or $\theta=3/4$.

The same issues arise when we test $\theta=1/2$ against the composite
hypothesis $\theta\ne1/2$. A natural Baye\-sian method for doing this is
to put half our probability on $\theta=1/2$ and distribute the other
half uniformly on $[0,1]$ (which is a special case of a widely
recommended procedure described in, e.g., \cite{bernardosmith2000},
page 391). This makes the likelihood ratio $X_t$ given
by (\ref{eq:process}) the posterior odds against $\theta=1/2$. As we
have seen, if the observations $x_1,x_2,\ldots$ turn out to be typical
for the distribution in which they are independent\vspace*{1pt} with the probability
for $x_t=1$ equal to $\frac12+\frac14\sqrt{\ln t/t}$, then if you wait
long enough, you can observe values of $X_t$ as small or as large as
you like, and thus obtain a posterior probability for $\theta=1/2$ as
large or as small as you like.

Of course, it will not always happen that the actual observations are
so equidistant from a simple null hypothesis and the probability
distribution representing its negation that the likelihood ratio will
oscillate wildly and you can sample to whichever side you want. More
often, the likelihood ratio and hence the posterior probability will
settle on one side or the other. But in the spirit of George Box's
maxim that all models are wrong, we can interpret this not as
confirmation of the side favored but only as confirmation that the
other side should be rejected. The rejection will be legitimate from
the Bayesian point of view, regardless of why we stopped sampling. It
will also be legitimate from the sampling-theory point of view.

On this argument, it is legitimate to collect data until a point has
been disproven but not legitimate to interpret this data as proof of an
alternative hypothesis within the model. Only when we really know the
model is correct can we prove one of its hypotheses by rejecting the
others.

\appendix
\section{Inadequacy of Test Martingales in Continuous Time}
\label{app:continuous}

In this appendix we will mainly discuss the case of continuous time; we
will see that in this case the notion of a test martingale is not fully
adequate for the purpose of hypothesis testing
(Proposition \ref{prop:local2}).
Fix a filtration $(\FFF_t)$ satisfying the usual conditions; in this
appendix we will only consider supermartingales $(X_t,\FFF_t)$, and we
will abbreviate $(X_t,\FFF_t)$ to $(X_t)$, or even to $X_t$ or~$X$.

In discrete time, there is no difference between\break using test martingales
and test supermartingales for hypothesis testing: every test martingale
is a test supermartingale, and every test supermartingale is dominated
by a test martingale (according to Doob's decomposition theorem,
\cite{meyer1966}, VII.1); therefore, using test supermartingales only
allows discarding evidence as compared to test martingales. In
continuous time, the difference between test martingales and test
supermartingales is essential, as we will see below
(Proposition~\ref{prop:local2}). For hypothesis testing we need ``local
martingales,'' a~modification of the notion of martingales introduced
by It\^{o} and Watanabe \cite{itowatanabe1965} and nowadays used
perhaps even more often than martingales themselves in continuous time.
This is the principal reason why in this article we use test
supermartingales so often starting from Section~\ref{sec:math}.

We will say that a random process $(X_t)$ is a \textit{local} member of a
class $\CCC$ of random processes (such as martingales or
supermartingales) if there exists a~sequence $\tau_1\le\tau_2\le\cdots$
of stopping times (called a~\textit{localizing sequence}) such that
$\tau_n\to\infty$ a.s. and each stopped process
$X^{\tau_n}_t=X_{t\wedge\tau_n}$ belongs to the class $\CCC$. (A
popular alternative definition requires that each
$X_{t\wedge\tau_n}\III_{\{\tau_n>0\}}$ should belong to $\CCC$.)
A standard argument (see, e.g.,
\cite{dellacheriemeyer1982}, VI.29) shows that there is no difference
between test supermartingales and local test supermartingales:
\begin{proposition}\label{prop:local1}
  Every local test supermartingale $(X_t)$
  is a test supermartingale.
\end{proposition}

\begin{pf}
  Let $\tau_1,\tau_2,\ldots$ be a localizing sequence,
  so that $\tau_n\to\infty$ as $n\to\infty$ a.s. and
  each $X^{\tau_n}$, $n=1,2,\ldots,$ is a test supermartingale.
  By Fatou's lemma for conditional expectations,
  we have, for $0\le s<t$,
  \begin{eqnarray*}
    \Expect(X_t\givn\FFF_s)
    &=&
    \Expect\Bigl(\lim_{n\to\infty} X^{\tau_n}_t\big\givn\FFF_s\Bigr)\\
    &\le&
    \liminf_{n\to\infty}
    \Expect(X^{\tau_n}_t\givn\FFF_s)\\
    &\le&
    \liminf_{n\to\infty}
    X^{\tau_n}_s
    =
    X_s\quad
    \mbox{a.s.}
  \end{eqnarray*}
  In particular, $\Expect(X_t)\le1$.
\end{pf}

An adapted process $(A_t)$ is called \textit{increasing} if $A_0=0$ a.s.\
and its every path is right-continuous
and increasing (as usual, not necessarily strictly increasing).
According to the Doob--Meyer decomposition theorem
(\cite{dellacheriemeyer1982}, Theorem VII.12), every test
supermartingale $(X_t)$ can be represented as the difference
$X_t=Y_t-A_t$ of a local test martingale $(Y_t)$ and an increasing
process $(A_t)$. Therefore, for the purpose of hypothesis testing in
continuous time, local test martingales are as powerful as test
supermartingales: every local test martingale is a test
supermartingale, and every test supermartingale is dominated by a local
test martingale. 

In discrete time there is no difference between local test martingales
and test martingales (\cite{dellacheriemeyer1982}, (VI.31.1)). In
continuous time, however, the difference is essential. Suppose the
filtration $(\FFF_t)$ admits a standard Brownian motion $(W_t,\FFF_t)$
in $\bbbr^3$. A well-known\vadjust{\goodbreak} example (\cite{johnsonhelms1963}; see also
\cite{meyer1966}, VI.21, and \cite{dellacheriemeyer1982}, VI.26) of a
local martingale which is not a martingale is
$L_t:=1/\|W_t+e\|$, where $e$ is a vector in $\bbbr^3$ such
that $\|e\|=1$ [e.g., $e=(1,0,0)$]; $L_t$ being a local
martingale can be deduced from $1/\|\cdot\|$ (the Newtonian
kernel) being a harmonic function on $\bbbr^3\setminus\{0\}$. The
random process $(L_t)$ is a local test martingale such that
$\sup_t\Expect(L_t^2)<\infty$; nevertheless, it fails to be a
martingale. See, for example, \cite{medvegyev2007} (Example 1.140) for
detailed calculations.

The local martingale $L_t:=1/\|W_t+e\|$ provides an example
of a test supermartingale which cannot be replaced, for the purpose of
hypothesis testing, by a~test martingale. According to another version
of the Doob--Meyer decomposition theorem (\cite{meyer1966}, VII.31), a
supermartingale $(X_t)$ can be represented as the difference
$X_t=Y_t-A_t$ of a martingale $(Y_t)$ and an increasing process $(A_t)$
if and only if $(X_t)$ belongs to the class (DL). The latter is defined
as follows: a supermartingale is said to be in (DL) if, for any $a>0$,
the system of random variables $X_{\tau}$, where $\tau$ ranges over the
stopping times satisfying $\tau\le a$, is uniformly integrable. It is
known that $(L_t)$, despite being uniformly integrable (as a collection
of random variables $L_t$), does not belong to the class (DL)
(\cite{meyer1966}, VI.21 and the note in VI.19). Therefore, $(L_t)$
cannot be represented as the difference $L_t=Y_t-A_t$ of a martingale
$(Y_t)$ and an increasing process $(A_t)$. Test martingales cannot
replace local test martingales in hypothesis testing also in the
stronger sense of the following proposition.

\begin{proposition}\label{prop:local2}
  Let $\delta>0$.
  It is not true that for every local test martingale $(X_t)$
  there exists a~test martingale $(Y_t)$
  such that $Y_t\ge\delta X_t$ a.s. for all~$t$.
\end{proposition}
\begin{pf}
  Let $X_t:=L_t=1/\|W_t+e\|$,
  and suppose there is a test martingale $(Y_t)$
  such that $Y_t\ge\delta X_t$ a.s. for all $t$.
  Let $\epsilon>0$ be arbitrarily small.
  Since $(Y_t)$ is in (DL)
  (\cite{meyer1966}, VI.19(a)),
  for any $a>0$ we can find $C>0$ such that
 \[
    \sup_{\tau}
    \int_{\{Y_{\tau}\ge C\}}
    Y_{\tau}\,
    \dd\Prob
    <
    \epsilon\delta,
 \]
  $\tau$ ranging over the stopping times satisfying $\tau\le a$.
  Since
  \[
    \sup_{\tau}
    \int_{\{X_{\tau}\ge C/\delta\}}
    X_{\tau}\,
    \dd\Prob
    \le
    \sup_{\tau}
    \int_{\{Y_{\tau}\ge C\}}
    (Y_{\tau}/\delta)\,
    \dd\Prob
    <
    \epsilon,
  \]
  $(X_t)$ is also in (DL),
  which we know to be false.
\end{pf}

\vspace*{-2pt}
\section{Details of Calculations}
\label{app:calculations}
\vspace*{-2pt}

\renewcommand{\theequation}{\arabic{equation}}
\setcounter{equation}{27}

In this appendix we will give details of some calculations omitted in
Section \ref{sec:example}. They will be based on Stirling's formula
$n!=\sqrt{2\pi n}(n/e)^n e^{\lambda_n}$, where $\lambda_n=o(1)$ as
$n\to\infty$.

\vspace*{-2pt}
\subsection{Oscillating Evidence when Testing Against a~Composite Alternative}
\vspace*{-2pt}

First we establish (\ref{eq:osc}) for $X_t$ defined by
(\ref{eq:process}). Suppose we have made $t$ observations and observed
$k:=k_t$ 1s so far. We start from finding bounds on $k$ that are
implied by the law of the iterated logarithm. Using the simplest
version of Euler's summation formula (as in \cite{apostol1999},
Theorem 1), we can find its expected value as
\begin{eqnarray*}
  \Expect(k)
  &=&
  \sum_{n=1}^t
  \Biggl( \frac12
    +
    \frac14
    \sqrt{\frac{\ln n}{n}} \Biggr)\\[-2pt]
  &=&
  \frac{t}{2}
  +
  \frac14
  \sum_{n=2}^t
  \biggl(
    \frac{\ln n+1}{\sqrt{n\ln n}}
  \biggr)
  -
  \frac14
  \sum_{n=2}^t
  \biggl(
    \frac{1}{\sqrt{n\ln n}}
  \biggr)\\[-2pt]
  &=&
  \frac{t}{2}
  +
  \frac14
  \int_{2}^t
  \biggl(
    \frac{\ln u+1}{\sqrt{u\ln u}}
  \biggr)\,
  \dd u
  +
  O\bigl(\sqrt{t}\bigr)\\[-2pt]
  &=&
  \frac{t}{2}
  +
  \frac12
  \sqrt{t\ln t}
  +
  O\bigl(\sqrt{t}\bigr).
\end{eqnarray*}
Its variance is
\begin{eqnarray*}
  \var(k)
  &=&
  \sum_{n=1}^t
  \Biggl(
    \frac12+\frac14\sqrt{\frac{\ln n}{n}}
  \Biggr)
  \Biggl(
    \frac12-\frac14\sqrt{\frac{\ln n}{n}}
  \Biggr)\\[-2pt]
  &=&
  \sum_{n=1}^t
  \biggl(
    \frac14-\frac{1}{16}\frac{\ln n}{n}
  \biggr)
  \sim
  \frac{t}{4}.
\end{eqnarray*}
Therefore, Kolmogorov's law of the iterated logarithm gives
\begin{eqnarray}\label{eq:k}
  \hspace*{25pt}\limsup_{t\to\infty}
  \frac
  {k-(1/2)(t+\sqrt{t\ln t})}
  {\sqrt{(1/2)t\ln\ln t}}
  &=&
  1\quad\mbox{and}\nonumber\\ [-9pt]\\ [-9pt]
  \liminf_{t\to\infty}
  \frac
  {k-(1/2)(t+\sqrt{t\ln t})}
  {\sqrt{(1/2)t\ln\ln t}}
  &=&
  -1\quad\mbox{a.s.}\nonumber
\end{eqnarray}

Using the definition (\ref{eq:process}) and applying Stirling's
formula, we obtain
\begin{eqnarray} \label{eq:by-Stirling}
  \ln X_t
  &=&
  t\ln2
  +
  \ln\frac{k!(t-k)!}{t!}
  -\ln(t+1)\nonumber\\
  &=&t\ln2-tH(k/t)
  +\ln\sqrt{2\pi\frac{k(t-k)}{t}}\nonumber\\
  &&{}+\lambda_k + \lambda_{t-k} - \lambda_t
  -\ln(t+1)\\
  &=&t\bigl(\ln 2 - H(k/t)\bigr)
  -\frac12\ln t+O(1)\nonumber\\
  &=&2t\biggl(\frac{k}{t}-\frac12\biggr)^2
  - \frac12\ln t+O(1)\quad\mbox{a.s.},\nonumber
\end{eqnarray}
where $H(p):=-p\ln p-(1-p)\ln(1-p)$, $p\in[0,1]$, is the entropy
function; the last equality in (\ref{eq:by-Stirling}) uses
$\ln2-H(p)=2(p-1/2)^2+O(|p-1/2|^3)$ as $p\to1/2$. Combining
(\ref{eq:by-Stirling}) with (\ref{eq:k}), we further obtain
%
\begin{eqnarray}\label{eq:cs}
  \limsup_{t\to\infty}
  \frac{\ln X_t}{\sqrt{2\ln t\ln\ln t}}
  &=&1\quad\mbox{and}\nonumber\\ [-8pt]\\ [-8pt]
  \liminf_{t\to\infty}
  \frac{\ln X_t}{\sqrt{2\ln t\ln\ln t}}
  &=&-1\quad\mbox{a.s.}\nonumber
\end{eqnarray}

%

\vspace*{-2pt}
\subsection{Prediction Interval}
\vspace*{-2pt}

Now we show that (\ref{eq:notrejm}) can be rewritten as
(\ref{eq:variable}). For brevity, we write $k$ for $k_n$. Similarly to
(\ref{eq:by-Stirling}), we can rewrite (\ref{eq:notrejm}) as
%
\begin{eqnarray}\label{eq:first}
  &&\hspace*{25pt}\ln2-H(k/n)+\frac1n
  \ln\sqrt{2\pi\frac{k(n-k)}{n}}
  \nonumber\\ [-8pt]\\ [-8pt]
  &&\hspace*{25pt}\quad{}+\frac{\lambda_k+\lambda_{n-k}-\lambda_n}{n}
  -\frac1n
  \ln(n+1)<\frac{\ln(1/\delta)}{n}.\nonumber
\end{eqnarray}
Since $\ln2-H(p)\sim2(p-1/2)^2$ ($p\to1/2$), we have $k/n=1/2+o(1)$ for
$k$ satisfying (\ref{eq:first}), as $n\to\infty$. Combining this with
(\ref{eq:first}), we further obtain
\begin{eqnarray*}
  &&2
  \biggl(
    \frac{k}{n} - \frac12
  \biggr)^2\\
  &&\quad<
  (1+\alpha_n)
  \frac{\ln(1/\delta)-\ln\sqrt{n}+\ln(n+1)+\beta_n}{n},
\end{eqnarray*}
for some $\alpha_n=o(1)$ and $\beta_n=O(1)$, which can be rewritten as
(\ref{eq:variable}) for a different sequence
$\alpha_n=o(1)$.

\subsection{Calculations for Armitage's Example}

Finally, we deduce (\ref{eq:cc-1}). Using a well-known expression
(\cite{abramowitzstegun1964}, 6.6.4) for the regularized beta function
$I_p(a,b):=B(p;a,b)/B(a,b)$ and writing $k$ for $k_t$, we obtain
\begin{eqnarray}\label{eq:cc-2}
X_t&=&\bigl(B(k+1,t-k+1)\nonumber\\
&&\phantom{\bigl(}{}-B(1/2;k+1,t-k+1)\bigr)\nonumber\\
&&{}/B(1/2;k+1,t-k+1)\nonumber\\ [-8pt]\\ [-8pt]
  &=&\frac{1}{I_{1/2}(k+1,t-k+1)}-1\nonumber\\
  &=&\frac{1}{\Prob\{B_{t+1}\ge k+1\}} -1\nonumber\\
 & =&\frac{\Prob\{B_{t+1}\le k\}}{\Prob\{B_{t+1}\ge k+1\}}.\nonumber
\end{eqnarray}

As a final remark, let us compare the sizes of oscillation of the log
likelihood ratio $\ln X_t$ that we have obtained in
Section \ref{sec:example} and in this appendix for our examples of the
three kinds of Bayesian hypothesis testing. When testing a simple null
hypothesis against a simple alternative, $\ln X_t$ oscillated between
approximately $\pm0.75\sqrt{t\ln\ln t}$ (as noticed in Section
\ref{subsec:ss}). When testing a simple null hypothesis against a
composite alternative, $\ln X_t$ oscillated between $\pm\sqrt{2\ln
t\ln\ln t}$ [see (\ref{eq:cs})]. And finally, when testing a composite
null hypothesis against a composite alternative, we can deduce from
(\ref{eq:cc-2}) that
\[
  \limsup_{t\to\infty}
  \frac{\ln X_t}{\ln\ln t}
  =
  1
  \quad\mbox{and}\quad
  \liminf_{t\to\infty}
  \frac{\ln X_t}{\ln\ln t}
  =
  -1\quad  \mbox{a.s.}
\]
(details omitted); therefore, $\ln X_t$ oscillates between $\pm\ln\ln
t$. Roughly, the size of oscillations of $\ln X_t$ goes down from
$\sqrt{t}$ to $\sqrt{\ln t}$ to $\ln\ln t$. Of course, these sizes are
only examples, but they illustrate a~general tendency.

\section*{Acknowledgments}

A. Philip Dawid and Steven de Rooij's help is gratefully appreciated.
Steven's thoughts on the subject of this article have been shaped by
discussions with Peter Gr\"{u}nwald. Comments by three reviewers have
led to numerous corrections and improvements, including addition of
Section \ref{sec:example}. We are grateful to Irina Shevtsova for
advising us on the latest developments related to the Berry--Esseen
theorem. In our computer simulations we have used the R language
\cite{R2010SS} and the GNU C++ compiler. Our work on the article has been
supported in part by ANR Grant NAFIT ANR-08-EMER-008-01 and EPSRC Grant
EP/F002998/1.

\end{document}